\documentclass[a4paper,11pt,twoside]{article}
\usepackage[latin1]{inputenc}
\usepackage[T1]{fontenc}
\usepackage[francais,english]{babel}
\usepackage{math,hog,graphicx}

\author{{\Large Vincent BEFFARA}\\
  Département de Mathématiques\\
  Université Paris-Sud\\
  F-91405 ORSAY Cedex, France\\
  \texttt{Vincent.Beffara@math.u-psud.fr}}
\date{}
\title{\Huge\bf On Conformally Invariant Subsets of the Planar Brownian Curve}

\renewcommand{\thmname}{Theorem}
\renewcommand{\dfname}{Definition}
\renewcommand{\lemname}{Lemma}
\hognewtheorem{cor}{Corollary}

\newcommand{\EA}[1][]{\ensuremath{\mathcal{E}_{A#1}}}
\newcommand{\EEA}[1][]{\ensuremath{\tilde{\mathcal{E}}_{A#1}}}
\newcommand{\TA}[1][]{\ensuremath{\mathcal{T}_{A#1}}}
\newcommand{\TTA}[1][]{\ensuremath{\tilde{\mathcal{T}}_{A#1}}}

\begin{document}

\maketitle 

\begin{abstract}
  We define and  study a family of generalized  non-intersection exponents for
  planar Brownian motions that is indexed by subsets of the complex plane: For
  each $A  \subset \CC$, we  define an exponent  $\xi (A)$ that  describes the
  decay  of   certain  non-intersection  probabilities.   To   each  of  these
  exponents,  we  associate  a  conformally  invariant subset  of  the  planar
  Brownian path,  of Hausdorff dimension $2-  \xi (A)$. A  consequence of this
  and continuity of $\xi(A)$ as a function of $A$ is the almost sure existence
  of  pivoting points of  any sufficiently  small angle  on a  planar Brownian
  path.
\end{abstract}

\selectlanguage{francais}

\begin{abstract}
  Nous  définissons et  étudions une  famille d'exposants  de non-intersection
  généralisés  entre mouvements browniens  plans, indexée  par les  parties du
  plan  complexe :  pour  chaque $A\subset\CC$  nous  définissons un  exposant
  $\xi(A)$ décrivant  la décroissance  de probabilités de  non-intersection. À
  chacun de ces exposants est associée une partie de la trajectoire brownienne
  qui est invariante sous l'action  des transformations conformes et qui a une
  dimension de Hausdorff égale à $2-\xi(A)$. Une conséquence de ce résultat et
  de la continuité  de $\xi(A)$ comme fonction de  $A$ est l'existence presque
  sûre  de points  pivotants de  tout angle  assez petit  sur  une trajectoire
  brownienne plane.
\end{abstract}

\selectlanguage{english}


\tableofcontents

\section*{Introduction}
\addcontentsline{toc}{section}{Introduction}
\markboth{INTRODUCTION}{INTRODUCTION}

It has been  conjectured for more than twenty  years by theoretical physicists
that conformal invariance plays an  important role to understand the behaviour
of critical two-dimensional models of  statistical physics.  They justify by a
mathematically non-rigorous  argument involving renormalization  ideas that in
the scaling  limit these  models behave in  a conformally invariant  way; they
have been able  to classify them via a  real-valued parameter corresponding to
the  central charge of  the associated  Virasoro algebra,  and to  predict the
exact  value  of critical  exponents  that  describe  the behaviour  of  these
systems.  Different models (for instance, self-avoiding walks and percolation)
with the same central charge have the same exponents.

Recently, Schramm \cite{schramm:UST} introduced  new mathematical objects that
give insight  into these conjectures.  These are  random set-valued increasing
processes  $(K_t)_{t\geqslant0}$ that  he called  Stochastic  Löwner Evolution
processes. For each positive number $\kappa$, there exists one such process of
parameter $\kappa$, in short $SLE_\kappa$.  He proved that for various models,
if they have a conformally invariant scaling limit, then it can be interpreted
in terms  of one of the  $SLE_\kappa$'s (the parameter $\kappa$  is related to
the central charge of the model).  One can then interpret the conjectures from
theoretical physicists in terms of properties of this process.

In particular,  Lawler, Schramm and  Werner  \cite{werner:value,werner:value2}
showed  that  for  one  specific  value  of  the  parameter  $\kappa$  (namely
$\kappa=6$)  which   conjecturally  corresponded  to  the   scaling  limit  of
percolation   cluster    interfaces,   the   $SLE_6$    has   the   remarkable
\emph{restriction property}  that enables to relate its  critical exponents to
the so-called  intersection exponents  between planar Brownian  motions.  This
lead \cite{werner:value, werner:value2,  werner:value3, werner:analyticity} to
the derivation  of the  exact value of  the exponents between  planar Brownian
paths.  Furthermore,  it turned out \cite{werner:barcelone} that  in fact, the
outer boundary of  a planar Brownian curve has exactly the  same law than that
of  an $SLE_6$.   In other  words,  the geometry  of critical  two-dimensional
percolation clusters in their scaling limit should be exactly that of a planar
Brownian outer frontier.

In a very recent  paper Smirnov \cite{smirnov:perco} showed that critical site
percolation in the triangular lattice  is conformally invariant in the scaling
limit so  that the geometry  of critical two-dimensional  percolation clusters
boundaries in their scaling limit is identical that of a planar Brownian outer
frontier.

Before all these recent  developments, geometric properties of planar Brownian
paths    had    already    been    subject   of    numerous    studies    (see
e.g. \cite{legall:sflour}  for  references).   In  particular,  the  Hausdorff
dimension of various subsets of the planar Brownian curve defined in geometric
terms had been determined.   For instance, Evans \cite{evans:cone} showed that
the Hausdorff dimension of the set  of two-sided cone points of angle $\theta$
(\ie\ points $B_t$ such that both $B_{[0,t]}$ and $B_{[t,1]}$ are contained in
the same  cone of angle $\theta$  with endpoint at $B_t$)  is $2-2\pi/\theta$. 
In a series of papers  (see \cite{lawler:fractal} for a review), Lawler proved
that the dimension  of various important subsets of  the planar Brownian curve
can     be    related    to     Brownian    intersection     exponents.     In
particular \cite{lawler:hausdorff}, he  showed that  the dimension of  the set
$C$ of cut  points (\ie\ points $B_t$ such  that $B_{[0,1]} \setminus \{B_t\}$
is  not connected) is  $2 -  \xi$ where  $\xi $  is the  Brownian intersection
exponent defined by
\begin{equation}
  \label{eq:isect}
  p_R = P(B^1_{[0,T^1_R]} \cap B^2_{[0,T^2_R]} = \varnothing) = R^{-\xi +o(1)}
\end{equation}
(for independent Brownian paths $B^1$ and $B^2$ starting respectively from $1$
and $-1$, $T^1_R$  and $T^2_R$ standing for their  respective hitting times of
the circle $\mathcal C(0,R)$).

In order  to derive such  results and in  particular the more  difficult lower
bound   $d  \geqslant   2-\xi$,  the   strategy   is  first   to  refine   the
estimate (\ref{eq:isect})  into $p_R  \asymp  R^{- \xi}$  (we  shall use  this
notation to denote  the existence of two positive constants  $c$ and $c'$ such
that $cR^{-\xi} \leqslant p_R  \leqslant c'R^{-\xi}$), to derive second-moment
estimates  and to  use these  facts to  construct a  random measure  of finite
$r$-energy  supported on  $C$, for  all $r<2-\xi$.   The determination  of the
value of the  critical exponents via $SLE_6$ \cite{werner:value,werner:value2}
then   implies   that   the   dimension   of   $C$   is   $3/4$.    Similarly,
in \cite{lawler:frontier} the  Hausdorff dimension of the outer  frontier of a
Brownian path  can be interpreted in  terms of another  critical exponent, and
the   determination   of   this    exponent   using   $SLE_6$   then   implied
(see \cite{werner:4/3}  for  a  review)   that  this  dimension  is  $4/3$  as
conjectured by Mandelbrot.
 
\bigskip

In the present  paper, we define and study a family  of generalizations of the
Brownian intersection  exponent $\xi$ parameterized by subsets  of the complex
plane. For each $ A \subset \CC$,  we define an exponent $\xi (A)$ as follows. 
Let $B^1$  and $B^2$  be two independent  planar Brownian paths  starting from
uniformly distributed points on the unit circle : then $\xi(A)$ is defined by
\begin{equation}
  \label{eq:defA}
  p_R(A) = P(B^1_{[0,T^1_R]} \cap A.B^2_{[0,T^2_R]} = \varnothing) =
  R^{-\xi(A) + o(1)}
\end{equation}
(with the notation $E_1.E_2 = \{xy:x\in E_1, y\in E_2\}$).  Note that the case
$A=\{1\}$    corresponds   to   the    usual   intersection    exponent.    In
Section \ref{sec:exponents}, we first show that for a wide class of sets $A$
\begin{equation}
  \label{eq:strong}
  p_R (A) \asymp R^{- \xi (A)}.
\end{equation}
  
In Section \ref{sub:propexp}, we study regularity properties of the mapping $A
\mapsto \xi (A)$.  In particular, we prove uniform continuity (with respect to
the Hausdorff  metric) on  certain families of  sets.  One important  tool for
this result is  the fact that the constants  implicit in (\ref{eq:strong}) can
in fact be taken uniform over these families of sets.
 
In Section \ref{sec:hausdorff}, we associate to each set $A$ a subset $\EA$ of
the planar Brownian curve defined in geometric terms:
$$\EA  =  \{B_t  :  \exists  \varepsilon>0,  (B_{[t-\varepsilon,t]}-B_t)  \cap
A.(B_{(t,t+\varepsilon]}-B_t)    =   \varnothing\}.$$
Using the  strong approximation and continuity  of the mapping  $A \mapsto \xi
(A)$, we then  show that the Hausdorff dimension of this  subset of the planar
Brownian curve is almost surely $2 - \xi  (A)$ (and $0$ in case $\xi (A) >2$). 
For example, when  $A=\{e^{i\theta}, 0\leqslant \theta\leqslant \alpha\}$, the
corresponding subset  $C_\alpha$ of the Brownian  curve is the  set of (local)
pivoting points, \ie\  points around which one half of the  path can rotate of
any angle smaller than $\alpha$ without intersecting the other half.

When $A \subset A'$, then  $\mathcal E_{A'} \subset \EA$.  In particular, when
$A$ contains $1$, then $\EA$ is a subset of the set of (local) cut points, and
therefore the shape of the path in a neighbourhood of such a point is the same
as the Brownian  frontier in the neighbourhood of  a cut-point.
This shows  in
particular that  (at least some of)  the exponents $\xi(A)$  describe also the
Hausdorff dimension  of sets of exceptional points 
of the scaling limit  of critical
percolation clusters.

In Section \ref{sec:numerical},  we derive some  bounds on the  exponents $\xi
(A)$  for  small  sets  $A$  in  the  same spirit  as  the  upper  bounds  for
disconnection exponents  derived in \cite{werner:bounds}.  In  particular, for
small $\alpha$, we show that  the exponent $\xi(C_\alpha)$ is strictly smaller
than $2$, which  implies the existence of pivoting points  (of small angle) on
the planar  Brownian curve.   We then briefly  present results  of simulations
that suggest that there exist pivoting points of angle up to an angle close to
$3  \pi /4$.  

\bigskip

Actually, it  is easy to define  other ``generalized'' exponents  in a similar
fashion, by studying non-intersection  properties between Brownian motions and
some of their images under isometries and scalings, \ie\ one can view $A$ as a
subset  of  the  linear   group.   Also,  one  can  consider  non-intersection
properties between $B$ and its image $f(B)$ by a conformal map.  It is easy to
see  using the  function  $z \mapsto  z^2$  that the  exponent describing  the
non-intersection  between $B$  and $-B$  is  in fact  twice the  disconnection
exponent.   The methods  of the  present  paper can  then be  adapted to  such
situations.

Similarly, one  could also  extend the definitions  to higher  dimensions (the
cases $d \geqslant  4$ can also be interesting if the  set $A$ is sufficiently
large), but conformal invariance can not  be used anymore, so that some of the
tools that we use in the present paper do not apply.

\subsection*{Acknowledgments}

I thank Wendelin \textsc{Werner} for suggesting me to look for pivoting points
on the planar Brownian curve and for never refusing help and advice.

\bigskip \bigskip

\subsection*{Notations}

Throughout this paper, we will  use the following notations for the asymptotic
behaviour of positive functions (and sequences, with the same meaning):
\begin{itemize}
\item $f\sim g$ if $\displaystyle \lim_{t\to\infty} \frac{f(t)}{g(t)} = 1$ ---
 and $f$ and $g$ are said to be \emph{equivalent};
\item   $f\approx  g$   if  $\log   f\sim\log  g$,   \ie\   if  $\displaystyle
  \lim_{t\to\infty} \frac{\log f(t)}{\log g(t)} =  1$ --- $f$ and $g$ are then
  \emph{logarithmically equivalent};
\item $f\asymp  g$ if  $f/g$ is bounded  and bounded  by below, \ie\  if there
  exist  two positive  finite constants  $c$ and  $C$ such  that for  all $t$,
  $cg(t)  \leqslant  f(t) \leqslant  Cg(t)$  ---  which  we call  \emph{strong
    approximation} of $f$ by $g$.
\end{itemize}


\section{Generalized intersection exponents}
\label{sec:exponents}

\subsection{Definition of the exponents}

\begin{hogtheorem}{Proposition and Definition :}
  Let $A$ be a  non-empty subset of the complex plane and  $B^1$, $B^2$ be two
  independent Brownian  paths starting uniformly on the  unit circle $\mathcal
  C(0,1)$; define the  hitting time $T_R^i$ of $\mathcal  C(0,R)$ by $B^i$ and
  let $\tau_n^i = T_{\exp(n)}^i$,
  $$E_n  =  E_n  (A)   =  \{B^1_{[0,\tau^1_n]}  \cap  A  B^2_{[0,\tau^2_n]}  =
  \varnothing \},$$  
  $$q_n (A) = P ( E_n) \hbox { and } p_R(A) = P(E_{\log R}).$$  
  Then, assuming  the existence  of positive constants  $c$ and $C$  such that
  $p_R(A) \geqslant cR^{-C}$,  there exists a real number  $\xi(A)$ such that,
  when $R\to\infty$,
  $$p_R(A) \approx R^{-\xi(A)}.$$
\end{hogtheorem}

\begin{dem}
  This is a standard sub-multiplicativity argument.  If $B$ is a Brownian path
  starting on  $\mathcal C(0,1)$  with \emph{any} law  $\mu$, then the  law of
  $B_{\tau_1(B)}$ on the  circle $\mathcal C(0,e)$ has a  density (relative to
  the  Lebesgue measure)  bounded  and  bounded away  from  zero by  universal
  constants  (\ie\ independently of  $\mu$).  Combining  this remark  with the
  Markov property  at the hitting  times of the  circle of radius  $e^n$ shows
  that:
  $$\forall m,n \geqslant 1 \qquad q_{m+n} \leqslant c q_n q_{m-1}.$$  
  Hence   the   family    $(cq_{n-1})$   is   sub-multiplicative   and   using
  Proposition \ref{prp:subadd}   we  have   $q_n\approx   e^{-\xi  n}$,   with
  $\xi\in(0,\infty)$, as well  as a lower bound $q_n  \geqslant c^{-2} e^{-\xi
    (n+1)}$.
\end{dem}

\emph{Remarks:}  For  some   choices  of  $A$  there  is   an  easy  geometric
interpretation  of   the  event   $E_n(A)$:  $\xi(\{1\})$  is   the  classical
intersection exponent; if  $A=(0,\infty)$, the $E_n(A)$ is the  event that the
paths stay in different wedges.

If $A$  is such that no lower  bound $p_R(A) \geqslant cR^{-C}$  holds, we let
$\xi(A)=\infty$.  However, in  most of  the  results presented  here, we  will
restrict ourselves to a class of sets  $A$ for which it is easy to derive such
lower bounds:

\begin{df}
  A non-empty subset $A$ of the complex  plane is said to be \emph{nice} if it
  is  contained   in  the  intersection  of  an  annulus  $\{r<|z|<R\}$  (with
  $0<r<R<\infty$) with a  wedge of angle strictly less  than $2\pi$ and vertex
  at $0$.
\end{df}

Indeed, let $A$  be such a set and  let $\alpha<2\pi$ be the angle  of a wedge
containing $A$: $B^1$ and $AB^2$ will not intersect provided each path remains
in a well-chosen wedge of angle  $(2\pi-\alpha)/2$, and then it is standard to
derive the following bound:
\begin{equation}
  \label{eq:apriori}
  p_R(A) \geqslant cR^{-4\pi/(2\pi-\alpha)}.
\end{equation}
The fact that $A$  be contained in an annulus will be  needed in the following
proof.  The only  usual case where this does  not hold is when $A$  is a wedge
itself; but in  this case a direct study is possible,  based on the derivation
of \emph{cone exponents} in \cite{evans:cone}  and the exact value of $\xi$ is
then known (cf. next section for details).

We will often consider the case where $A$ is a subset of the unit circle.  For
such sets, $A$  is nice if and only  if $\bar A \varsubsetneq \dUU$  (it is in
fact easy  to prove  that for $A\subset\dUU$,  $\xi(A)=\infty$ if and  only if
$\bar A=\dUU$).

\subsection{Strong approximation}

This  whole subsection  will be  dedicated  to the  refinement of  $p_R\approx
R^{-\xi}$  into $p_R\asymp  R^{-\xi}$. This  is \emph{not}  anecdotical, since
this ``strong'' approximation  will be needed on several  occasions later.

\begin{thm}
  \label{thm:strong}
  For every  nice $A$, $p_R(A)\asymp  R^{-\xi(A)}$, \ie\ there  exist positive
  constants $c(A)<C(A)$ such that
  $$c R^{-\xi(A)} \leqslant p_R(A) \leqslant CR^{-\xi(A)}.$$  
  Moreover,  the constants  $c(A)$  and $C(A)$  can  be taken  uniformly on  a
  collection $\mathcal  A$ of subsets of  the plane, provided  the elements of
  $\mathcal A$ are contained in the same nice set.
\end{thm}

\begin{dem}
  Note that since $A \in {\cal A}$ is nice, the exponents $\xi (A)$ exists and
  is uniformly  bounded (for  $A \in {\cal  A}$).  The  subadditivity argument
  showed that $q_n \geqslant c  e^{-\xi(A). (n+1)}$, which implies readily the
  lower bound in the theorem.  It is more difficult to derive the upper bound.
  By  Proposition \ref{prp:subadd}, it  will be  sufficient to  find  a finite
  constant $c_-(A)$ (that can be bounded  uniformly for $A \in {\cal A}$) such
  that
  \begin{equation}
    \label{eq:need}
    \forall n,n' \qquad q_{n+n'} \geqslant c_- q_n q_{n'}. 
  \end{equation}
  In order to make the proof more readable, it is carried out here for a fixed
  $A$; however  it is easy  to see  that, at each  step, the constants  can be
  taken  uniformly  for all  $A$  contained  in some  fixed  nice  set $A\z$.  
  Moreover, we shall  first assume that $A\z$ is a subset  of the unit circle:
  We briefly indicate  at the end of the proof what  are the few modifications
  needed to adapt it to the general case.
  
  The  basic  method  is  adapted  from Lawler's  proof  for  non-intersection
  exponents  in \cite{lawler:concavity},  with some  technical simplifications
  made possible using the absence of the $\lambda$ exponent.  The main idea is
  to obtain a weak independence between  the behaviour of the paths before and
  after they reach radius $e^n$.  The first step is an estimate concerning the
  probability that  the paths  are ``well separated''  when they  reach radius
  $e^n$  (more precisely,  that they  remain in  two separated  wedges between
  radius $e^{n-1}$ and radius $e^n$):

  \begin{inlem}[Technical]
    Let $\eta>0$  and $\alpha  < 2\pi-\eta$  such that $A$  is contained  in a
    wedge of angle less than $\alpha$. Define 
    $$W_\alpha = \left\{ re^{i\theta} : r>0, |\theta|<\frac\alpha2 \right\},$$
    $\delta_n=    e^{-n}   [\dd(B^1_{\tau^1_n},    AB^2_{[0,\tau^2_n]})   \vee
    \dd(AB^2_{\tau^2_n}, B^1_{[0,\tau^1_n]})]$ and the following events:
    \begin{gather*}
      U^1_n = \left\{ B^1_{[0,\tau^1_n]} \cap \{|z|\geqslant e^{n-1}\} \subset
        -W_{2\pi-\alpha-\eta} \right\}, \\
      U^2_n  =  \left\{   AB^2_{[0,\tau^2_n]}  \cap  \{|z|\geqslant  e^{n-1}\}
        \subset W_\alpha \right\},
    \end{gather*}
    and $U_n=U^1_n \cap U^2_n$. Then:
    $$\exists  c,\beta>0\quad \forall\varepsilon>0\quad  \forall r  \in \left[
      \frac32,    3\right]\qquad    P(E_{n+r},U_{n+r}|E_n,\delta_n   \geqslant
    \varepsilon) \geqslant c\varepsilon^\beta.$$
  \end{inlem}
  
  \begin{dem}
    This is  an direct consequence of classical  estimates concerning Brownian
    motion in wedges; the value of  $\beta$ is not important, so not much care
    is needed in finding the lower  bound. Note that the existence of $\alpha$
    requires that $A$ be nice.
  \end{dem}
  
  If $\mathcal F_n$  stands for the $\sigma$-field generated  by both paths up
  to radius  $e^n$ (so that for instance  $E_n$ is in $\mathcal  F_n$), we now
  prove that paths conditioned not to  intersect up to radius $e^{n+2}$ have a
  good chance to  be well separated at this radius,  uniformly with respect to
  their behaviour up to radius $e^n$:

  \begin{inlem}[End-separation]
    There exists $c>0$ such that, for every $n>0$:
    $$P(U_{n+2}|E_{n+2},\mathcal F_n) \geqslant c$$     
    (\ie\ the  essential lower bound of  $P(U_{n+2}|E_{n+2},\mathcal F_n)$, as
    an $\mathcal F_n$-measurable function, is not less than $c$).
  \end{inlem}

  \begin{dem}
    The technical  lemma states that  start-separation occurs if  the starting
    points are sufficiently  far from each other; more  precisely, we have for
    all $\varepsilon>0$:
    \begin{equation}
      \label{eq:tech}
      P(U_{n+2}|E_{n+2},\mathcal F_n,\delta_n\geqslant\varepsilon) \geqslant
      c \varepsilon^\beta.
    \end{equation}
    Hence, what is to be proved is that two paths conditioned not to intersect
    have a positive  probability to be far from each  other after a relatively
    short time. To  prove this fact, one has to use  conditioning on the value
    of $\delta_n$.
    
    Fix $k>0$, and  assume that $2^{-(k+1)} \leqslant \delta_n  < 2^{-k}$; let
    $\tau_k$  be the  smallest $r$  such that  one of  the following  happens:
    either $\delta_{n+r} \geqslant 2^{-k}$, or  $E_{n+r}$ does not hold. It is
    easy to use scaling to prove that for some $\lambda>0$,
    $$P(\tau_k \geqslant 2^{-k}) \leqslant 2^{-\lambda},$$     
    meaning that  with positive probability  (independent of $k$ and  $n$) the
    paths separate  or meet before  reaching radius $e^{n+2^{-k}}$.   Hence by
    the strong Markov property, applying this $k^2$ times leads to
    \begin{equation}
      \label{eq:sepk}
      P(\tau_k \geqslant k^2 2^{-k}) \leqslant 2^{-\lambda k^2}.
    \end{equation}
    The technical lemma  states that $P(E_{n+2}|\delta_n \geqslant 2^{-(k+1)})
    \geqslant c 2^{-\beta k}$: combining both estimates then leads to
    \begin{equation}
      \label{eq:sepk2}
      P(\tau_k \geqslant k^2 2^{-k} | E_{n+2}, \delta_n \geqslant 2^{-(k+1)})
      \leqslant c 2^{\beta k - \lambda k^2}.
    \end{equation}
    
    Consider now a generic  starting configuration at radius $e^n$, satisfying
    $E_n$ and  hence $\delta_n>0$.  Fix also $k_0>0$  and introduce  the radii
    $\tau_k$ (for $k_0 \leqslant k < \infty$) defined by
    $$\tau_k = \Inf \{ r : \delta_{n+r} \geqslant 2^{-k} \}$$        
    (so   that   $\tau_k=0$   as    long   as   $2^{-k}\leqslant   \delta$).   
    Equation~(\ref{eq:sepk2})  can  be  rewritten  (using the  fact  that  the
    technical lemma is valid for all $r\geqslant3/2$) as
    $$P(\tau_{k}    -   \tau_{k+1}   \geqslant    k^2   2^{-k}    |   E_{n+2},
    \tau_{k+1}\leqslant \frac12) \leqslant c 2^{\beta k - \lambda k^2}.$$
    Fix $k_0$ such that
    $$\sum_{k=k_0}^\infty k^2 2^{-k} < \frac12,$$
    and sum this estimate for $k_0 \leqslant k < \infty$: this leads to
    $$P(\forall  k\geqslant  k_0,  \tau_k-\tau_{k+1}  \leqslant k^2  2^{-k}  |
    E_{n+2})  \geqslant 1  - c  \sum_{k=k_0}  ^{\infty} 2^{\beta  k -  \lambda
      k^2}.$$    
    In particular, if $k_0$ is taken large enough, this probability is greater
    than $1/2$, and we obtain
    $$P(\tau_{k_0} \leqslant \frac12 | E_{n+2}) \geqslant \frac12.$$
    It is then sufficient to combine this and Equation~(\ref{eq:tech}) to get
    $$P(U_{n+2} | E_{n+2}) \geqslant c 2^{-\beta k_0} >0,$$    
    and  is can be  seen that  the obtained  constant does  not depend  on the
    configuration at radius $e^n$ --- provided $E_n$ is satisfied.
  \end{dem}
  
  The first consequence of the end-separation lemma is $P(E_n,U_n)\asymp q_n$;
  but it  is easy to see, using  estimates on Brownian motion  in wedges again
  and the strong Markov property, that
  $$P(E_{n+1}|E_n,U_n)\geqslant c>0$$  
  (with  $c$  independent of  $n$),  and  combining  both estimates  leads  to
  $q_{n+1} \geqslant cq_n$, \ie\ $q_{n+1}\asymp q_n$. Now if $\bar q_n$ stands
  for the upper bound for the non-intersection probabilities, namely
  $$\bar q_n \eqd \Sup_{B^1\z, B^2\z \in \UU} P(E_n|B^1\z, B^2\z),$$  
  the previous  remark concerning  the law of  $W_{\tau_1(W)}$ can be  used to
  prove that $\bar q_n \leqslant cq_{n-1}$: hence,
  $$\bar q_n \asymp q_n.$$
  
  \bigskip
  
  Now that we know that paths  conditioned not to intersect have a good chance
  to exit  a disk  at a  large distance from  each other,  what remains  to be
  proven is  that paths  starting from distant  points on  $\mathcal C(0,e^n)$
  remain well separated for a sufficiently long time and become (in a sense to
  be specified later) independent from their behaviour before radius $e^n$.

  \begin{inlem}[Start-separation]
    Let $\alpha$ and  $\eta$ be as in the  technical lemma, $\eta'=\eta/2$ and
    $\alpha' = (2\pi + \alpha)/2$; introduce
    \begin{gather*}
      J^1_n = \left\{ B^1_{[0,\tau^1_n]} \cap \mathcal B(0,2) \subset
        -W_{2\pi-\alpha'-\eta'} \setminus \mathcal B(0,1-\eta') \right\}, \\
      J^2_n  =  \left\{   AB^2_{[0,\tau^2_n]}  \cap  \mathcal  B(0,2)  \subset
        W_{\alpha'} \setminus \mathcal B(0,1-\eta') \right\},
    \end{gather*}
    and $\tilde E_n=E_n\cap J^1_n \cap J^2_n$. Define $\tilde q_n$ as
    $$\tilde q_n(x,y) = P(\tilde E_n|B^1\z=x, B^2\z=y).$$    
    Then there exists $c>0$ such  that, for all $n\geqslant2$ and uniformly on
    all  pairs $(x,y)$  satisfying  $U\z$  (\ie\ such  that  $U\z$ holds  when
    $B^1\z=x$ and $B^2\z=y$):
    $$\tilde q_n(x,y) \geqslant cq_n.$$
  \end{inlem}
  
  \begin{dem}
    Introduce the following (``forbidden'') sets:
    \begin{align*}
      K^1 &= \left( \mathcal B(0,e) \setminus -W_{2\pi-\alpha'-\eta'} \right)
      \cup  \mathcal B(0,1-\eta');\\
      K^2 &= \left( \mathcal B(0,e) \setminus W_{\alpha'} \right) \cup
      \mathcal B(0,1-\eta').
    \end{align*}
    For all $n$ we have $J^1_n  = \{ B^1_{[0,\tau^1_n]} \cap K^1 = \varnothing
    \}$ and  $J^2_n = \{ AB^2_{[0,\tau^2_n]}  \cap K^2 =  \varnothing \}$. For
    the rest of the  proof we shall fix $n$, and condition  the paths by their
    starting  points; introduce  the  following stopping  times (for  positive
    values of $k$):
    \begin{align*}
      T^1\z &= \Inf \{ t>0 : B^1_{[0,t]} \cap \mathcal C(0,3) \neq
      \varnothing\},\\
      S^1_k &= \Inf \{ t>T^1_{k-1} : B^1_{[T^1_{k-1},t]} \cap K^1 \neq
      \varnothing \},\\
      T^1_k &= \Inf \{ t>S^1_k : B^1_{[S^1_k,t]} \cap \mathcal C(0,3) \neq
      \varnothing \},
    \end{align*}
    and  $S^2_k$, $T^2_k$  similarly, replacing  all occurrences  of  $B^1$ by
    $AB^2$ and $K^1$  by $K^2$.  We shall also use the  notation $N^i$ for the
    \emph{number  of crossings}  by $B^1$  (resp.\ $AB^2$)  between  $K^i$ and
    $\mathcal C(0,3)$, defined as
    $$N^i = \Max \{ k : S^i_k  < \tau^i_n \}.$$    
    With those notations, $J^i_n = J^i_1 \cap \{ N^i_n = 0\}$ and a.s.\ $N^i <
    \infty$.  Moreover,  uniformly  on  the starting  points  considered  here
    (satisfying the condition  $U\z$), we have $J^i_1 \geqslant c  > 0$ by the
    technical lemma, where $c$ depends only on $\eta$. 
    
    First, we split the event $E_n$  according to the value of, say, $N^2$: we
    write $P(E_n) = \sum_{k=0}^\infty P(E_n,N^2=k)$. By the Beurling estimate,
    on  $\{N\geqslant  k\}$,  the  probability that  $B^1_{[0,\tau^1_n]}$  and
    $AB^2_{[S^2_k,T^2_k]}$  do  not intersect  is  bounded  by some  universal
    constant  $\lambda<1$  (which can  even  be  chosen  independent of  $A$),
    independently  of $B^1$  and the  two remaining  parts of  $B^2$.   By the
    strong  Markov property at  time $T^2_{k}$,  when $N^2=k$  the probability
    that $AB^2$ after $T^2_{k}$ does  not intersect $B^1$ is bounded by $P(B^1
    \cap AB^2_{[T^2\z,\tau^2_n]} = \varnothing  , N^2=0)$ (\ie\ the path after
    $T^2_{k}$  when $N^2=k$  is the  same as  the entire  path when  $N^2=0$). 
    Introducing those two estimate in the sum leads to
    $$P(E_n) \leqslant  \sum_{k=0}^{\infty} \lambda^k P(E_n, N^2=0)  = \frac 1
    {1-\lambda} P(E_n,  N^2=0).$$    
    Doing this decomposition again according  to $N^1$ (with the same constant
    $\lambda<1$) we then obtain
    $$P(E_n) \leqslant \frac 1 {(1-\lambda)^2} P(E_n,N^1=N^2=0),$$    
    \ie\  $P(N^1=N^2=0|E_n)  \geqslant  (1-\lambda)^2  > 0$.   This,  and  the
    previous remark that $P(J^i_n | N^i=0)$  is bounded by below by a constant
    provided that the  starting points satisfy $U\z$, gives:
    \begin{equation}
      \label{eq:presque}
      P(\tilde E_n|B^1\z=x, B^2\z=y) \geqslant c P(E_n|B^1\z=x, B^2\z=y).
    \end{equation}

    Conditioning on $B^2$ shows that the map
    \begin{equation}
      \label{eq:func}
      f:x \mapsto P(E_n|B^1\z=x, B^2\z=1)
    \end{equation}
    is  harmonic and  does not  vanish  on the  complement of  $\overline A$.  
    Moreover,  its supremum  on the  unit  circle is  equal to  $\bar q_n$  by
    definition: Applying the Harnack principle then proves that $f$ is bounded
    by below by $cq_n$ on the set of $x$ satisfying $U\z$, which completes the
    proof.
  \end{dem}
  
  Another estimate can be obtained using the very same proof: Only keeping the
  conditions involving disks and relaxing those involving wedges, we obtain
  \begin{equation}
    \label{eq:nowedge}
    P \left( B^1_{[0,\tau^1_n]} \cap B(0,1-\eta) = \varnothing,
      AB^2_{[0,\tau^2_n]}  \cap \mathcal B(0,1-\eta) = \varnothing \Big|
      B^1\z, B^2\z, E_n \right) \geqslant c>0,
  \end{equation}
  where $c$ does not depend on  the initial positions $B^1\z$ and $B^2\z$, nor
  on $n$ (it clearly depends on $\eta$, though, and a closer look at the proof
  shows that we can ensure  $c>\eta^\beta$ as $\eta\to0$, for some $\beta>0$). 
  This  estimate will  be needed  in the  derivation of  Hausdorff dimensions,
  cf. Section \ref{sec:hausdorff}.

  \bigskip
  
  We  know have all  the needed  estimates to  derive the  lower bound  in the
  sub-additivity condition, and hence the  conclusion of the theorem. Take two
  paths  with independent starting  points uniformly  distributed on  the unit
  circle and killed at radius  $e^{m+n}$, conditioned not to intersect between
  radii  $1$ and  $e^n$.  This  happens  with probability  $q_n$.  With  large
  probability (\ie\ with  a positive probability, independent of  $m$ and $n$)
  the paths up to  radius $e^n$ end up ``well separated'' in  the sense of the
  end-separation  lemma.  In particular,  the points  where they  reach radius
  $e^n$,   after   suitable  rescaling,   satisfy   the   hypothesis  of   the
  start-separation  lemma: Hence  with  probability greater  that $cq_m$,  the
  paths  between radii  $e^n$  and  $e^{m+n}$ remain  separated  up to  radius
  $e^{n+1}$, do not reach radius $(1-\eta)e^n$ anymore and do not intersect up
  to radius  $e^{m+n}$.  Under those  conditions, it is  easy to see  that the
  paths  do not  meet  at all.   So $q_{m+n}  \geqslant  c q_m  q_n$ for  some
  positive $c$, and we get the conclusion.
  
  \bigskip
  
  Some  adaptations  are  needed  if  $A$  is  included  in  an  annulus,  say
  $\{r<|z|<R\}$ with $r<1<R$. First, replace  all occurrences of $e$ by $e\z$,
  with $e\z$  chosen larger than  $10R/r$, and in the  start-separation lemma,
  replace $\mathcal B(0,1-\eta)$ by  $\mathcal B(0,r/2R)$ in the definition of
  the $J_n$.  As long as $r$  and $R$ are  fixed, this changes nothing  to the
  proof, except  that the constants  we obtain will  then depend on  $R/r$ ---
  which itself  is bounded provided  $A$ remains a  subset of some  fixed nice
  set. 
  
  A more  serious problem  arises if  the complement of  $\overline A$  is not
  connected,  since the  natural domain  of the  function $f$  (as  defined by
  Equation (\ref{eq:func})) is  itself not  connected.  However, since  $A$ is
  nice, its complement has exactly one  unbounded component, and it is easy to
  see  that  if  $x$  is  not  in this  component  then  $f(x)$  vanishes  for
  $n\geqslant  1$.   Hence,  nothing   changes  (as  far  as  non-intersection
  properties  are concerned) when  $A$ is  replaced by  the complement  of the
  infinite component  of its  complement (\ie\ when  filling the  ``holes'' in
  $A$).
\end{dem}

In fact, a stronger result can  be derived: If the starting points $B^1\z$ and
$B^2\z$   are  fixed,  then   $P(E_n|B^1\z,B^2\z)$  is   \emph{equivalent}  to
$ce^{-n\xi(A)}$, where  $c$ is  a function of  $B^1\z$ and  $B^2\z$ satisfying
$c\leqslant c\z  d(B^1\z,AB^2\z)^\beta$. This estimate is related  to a strong
convergence   result   on  the   law   of   paths   conditioned  by   $B^1\cap
AB^2=\varnothing$. However,  proving this result  would be much  more involved
(cf. \cite{lawler:intersection} for the proof in the case $A=\{1\}$).


\section{Properties of the function $A\mapsto\xi(A)$}
\label{sub:propexp}

We first  list a few simple  properties of the  function $A \mapsto \xi  (A)$. 
For  $p\in\ZZ$  and  $A\subset\CC$,  introduce $A^p=\{z^p,z\in  A\}$  and  let
$A^*=\{\bar z, z\in A\}$.

\begin{prp}
  \label{prp:props}
  Is these statements, all sets are assumed to be non-empty but do not need to
  be nice:
  \begin{enumerate}
  \item \vspace{-\lastskip}
  $\xi$  is  \emph{non-decreasing}:   if  $A\subset  A'$  then  $\xi(A)
    \leqslant \xi(A')$;
  \item $\xi$ is \emph{homogeneous}: if $\lambda\in\CC^*$ then $\xi(\lambda A)
    = \xi(A)$;
  \item $\xi$ is \emph{symmetric}: $\xi(A^{-1})=\xi(A^*)=\xi(A);$
  \item $\xi$ has the following property: if $n\geqslant 1$ then
    $$\xi\left(\bigcup e^{2ik\pi/n} A\right) = n\xi(A^n).$$
  \end{enumerate}
  \vspace{-\lastskip}
\end{prp}

\begin{dem}
  (i): This is a trivial consequence of $p_R(A) \geqslant p_R(A')$.
  
  (ii): Applying the scaling property  with factor $|\lambda|$ to $B^2$ proves
  that one  can suppose  $|\lambda|=1$; in  which case we  have $p_R(A)  = p_R
  (\lambda A)$ (because  the starting points are uniformly  distributed on the
  unit circle).
  
  (iii):  Simply exchange  $B^1$  and $B^2$  for  $A^{-1}$, and  say that  the
  complex conjugate of a Brownian path is still a Brownian path to get $A^*$.
  
  (iv): This is a consequence of the analyticity of the mapping $z\mapsto z^n$
  (hence the fact that $((W_t)^n)$ is  a Brownian path if $W$ is one) together
  with  the remark  that  the existence  of  $s$, $t>0$  and  $z\in A^n$  with
  $(B^1_s)^n = z (B^2_t)^n$ is equivalent to the existence of $z'$ in $\bigcup
  e^{2ik\pi/n} A$ with  $B^1_t=z'B^2_t$ --- note that the  mapping also has an
  influence on $R$, hence the factor $n$.
\end{dem}

We  now turn our  attention toward  regularity properties  of the  function $A
\mapsto  \xi(A)$  ---  the following  result  being  a  key step  towards  the
derivation  of  dimensions  in  the  next section.   Introduce  the  Hausdorff
distance       between      compact       subsets      of       the      plane
(cf. Section \ref{sec:ingredients} for  details).  It will  be convenient here
to define neighbourhoods by $V_r(A)=\{xe^{z},  x\in A, |z|<r\}$ instead of the
usual  $A+\mathcal B(0,r)$  ---  leading to  the \emph{logarithmic}  Hausdorff
distance.  The (logarithmic) Hausdorff topology is the metric topology derived
from this distance.

\begin{prp}
  \label{prp:cont}
  $\xi$  is  continuous on  the  collection of  nice  sets,  endowed with  the
  logarithmical  Hausdorff  topology.   For  any  nice  set  $A_0$,  $\xi$  is
  uniformly continuous in $\{ A \ : \ A \subset A_0 \}$.
 \end{prp}

\begin{dem}
  The  proof  relies  on  the   uniformity  of  the  strong  approximation  in
  Theorem \ref{thm:strong}:  fix  a  nice   set  $A\z$  and  assume  all  sets
  considered here  are subsets  of $A\z$. The  constants $c$, $c_-$  and $c_+$
  appearing during the proof may only depend on $A\z$.
  
  First, fix  $R>1$ and condition all  events by \smash{$B^2_{[0,T^2_{R+1}]}$}
  --- \ie\ fix the second path.  For all $A\subset A\z$, let
  $$\dd_R(A) = \dd(B^1 _{[0,T^1_R]},  AB^2 _{[0,T^2_R]})\;;$$
  for all  $\varepsilon>0$  introduce  the  stopping  time
  $$S_\varepsilon =  \Inf\{t: \dd(B^1_t,AB^2_{[0,T^2_R]}) <  \varepsilon \}.$$
  Note that $\{\dd_R(A)<\varepsilon\} = \{ S_\varepsilon < T^1_R \}$.  On this
  event, the strong Markov property shows that $B^1_{S_{\varepsilon} + \cdot}$
  is a  Brownian path starting  $\varepsilon$-close to $AB^2$.   By Beurling's
  theorem,  the probability  that  they do  not  meet before  radius $R+1$  is
  smaller  than the corresponding  probability for  a path  near a  half line;
  hence,
  $$P(B^1_{[S_\varepsilon,T^1_{R+1}]}   \cap    A   B^2_{[0,   T^2_{R+1}]}   =
  \varnothing | \dd_R(A) <  \varepsilon ) \leqslant \sqrt{\varepsilon},$$  
  so  that, considering  the  whole path,  $P(E_{R+1}|\dd_R(A) <  \varepsilon)
  \leqslant \sqrt{\varepsilon}$. Apply the Bayes formula:
  $$P(\dd_R(A)  <  \varepsilon |  E_{R+1})  = \frac  {P(\dd_R(A)<\varepsilon)}
  {P(E_{R+1})} P(E_{R+1}|\dd_R(A)  < \varepsilon);$$    
  since  we know  that $P(E_{R+1})\geqslant  c_-(R+1)^{-\xi(A)}$  with $\xi(A)
  \leqslant \xi(A\z)$ we finally obtain
  $$P(\dd_R(A)   <   \varepsilon   |   E_{R+1})   \leqslant   c   R^{\xi(A\z)}
  \sqrt{\varepsilon}.$$
  
  From now  on, we  shall assume that  $\varepsilon$ is sufficiently  small to
  make the obtained bound smaller that $1$. Taking the complement leads to
  $$P(\dd_R(A) \geqslant  \varepsilon | E_{R+1}) \geqslant 1  - c R^{\xi(A\z)}
  \sqrt{\varepsilon}.$$    
  Now, remark  that when $\dd_R(A)  \geqslant \varepsilon$ and  $\dd_H(A,A') <
  \varepsilon/R$,   we   have   $B^1_{[0,T^1_R]}  \cap   A'B^2_{[0,T^2_R]}   =
  \varnothing$: from this  and the previous equation follows  that, as long as
  $A$ and $A'$ remain subsets of $A\z$,
  $$\dd_H(A,A')  <  \frac  \varepsilon  R \quad\imply\quad  p_R(A')  \geqslant
  \left( 1-c R^{\xi(A\z)} \sqrt{\varepsilon} \right) p_{R+1}(A).$$
  
  We can apply  the estimates on $p_R$ we  derived in Theorem \ref{thm:strong}
  --- \ie\   $p_R(A)\asymp   p_{R+1}(A)   \asymp   R^{-\xi(A)}$:   still   for
  $\dd_H(A,A')<\varepsilon/R$ and $A$, $A'$ inside $A\z$ we get
  $$c_+  R^{-\xi(A')}  \geqslant  \left( 1-c  R^{\xi(A\z)}  \sqrt{\varepsilon}
  \right) c_- R^{-\xi(A)},$$  
  and taking the logarithm of each side of the inequality leads to
  $$\log  c_+  -  \xi(A')  \log  R  \geqslant  \log  c_-  +  \log  \left(  1-c
    R^{\xi(A\z)}  \sqrt{\varepsilon} \right)  - \xi(A)  \log R,$$  
  hence after suitable transformations:
  \begin{equation}
    \label{eq:unif}
    \xi(A')  \leqslant \xi(A)  + \frac  c {\log  R} -  \frac {\log  \left( 1-c
      R^{\xi(A\z)} \sqrt{\varepsilon}\right)} {\log R}.
  \end{equation}
  
  Fix  $\eta>0$, and choose  $R$ such  that $c/\log  R <  \eta/2$. It  is then
  possible  to take $\varepsilon$  sufficiently small  so that  $|\log (1  - c
  R^{\xi(A\z)}  \sqrt{\varepsilon})|  < (\eta\log  R)/2$;  for $\dd_H(A,A')  <
  \varepsilon/R$  we then  have $\xi(A')  \leqslant \xi(A)  + \eta$,  hence by
  symmetry  $|\xi(A') -  \xi(A)| \leqslant  \eta$. This  proves that  $\xi$ is
  uniformly continuous on $\mathcal P_c(A\z)$, for all $A\z$, hence continuous
  on the family of nice sets.
\end{dem}

\emph{Remark 1:} Equation (\ref{eq:unif}) allows the derivation of an explicit
modulus of  continuity for  $\xi$ inside $A\z$, of the form
$$|\xi(A')  -  \xi(A)| \leqslant  \frac  {C(A\z)} {|\log\dd_H(A,A')|}$$
(take $R=d^{-1/2\xi(A\z)}$).  But since $C(A\z)$  is not known, this  does not
provide numerical bounds for $\xi$.

\emph{Remark  2:} Inside  a  nice  set, the  usual  and logarithmic  Hausdorff
topologies   are   equivalent,    so   the   introduction   of   ``exponential
neighbourhoods'' in  Proposition \ref{prp:cont} can seem  artificial; however,
it leads to constants that do not vary when $A$ is multiplied by some constant
(as  in  Proposition \ref{prp:props}, point  (ii)),  hence uniform  continuity
holds on the  collection of nice sets  contained in a fixed wedge  and in some
annulus  $\{r<|z|<cr\}$  for  fixed $c$  ---  which  is  wrong for  the  usual
Hausdorff topology, as  a consequence of the homogeneity  of $\xi$ applied for
small $|\lambda|$.

Note that uniform continuity cannot hold  on the family of nice sets contained
in  a  given  annulus since  $\xi$  would  then  be  bounded (by  a  compacity
argument), which it is not: the exponent associated to a circle is infinite.


\section{Hausdorff dimension of the corresponding subsets of the path}
\label{sec:hausdorff}

\subsection{Conformally invariant subsets of the Brownian path}

It is  well-known that the Brownian  path is invariant in  law under conformal
transformations; in this section, we  study subsets of the Brownian curve that
are  also invariant  under conformal  maps.   A first  example is  the set  of
so-called Brownian \emph{cut-points}, \ie\  points $B_t$ such that $B_{[0,t)}$
and $B_{(t,1]}$ are  disjoint; these points form a  set of Hausdorff dimension
$2-\xi(\{1\})=3/4$. Related to those  are \emph{local cut-points}, \ie\ points
such that  there exists $\varepsilon>0$  satisfying $B_{[t-\varepsilon,t)}\cap
B_{(t,t+\varepsilon]}  = \varnothing$  --- the  dimension is  the same  as for
global     cut-points.     Other    examples     are    given     by    Lawler
in \cite{lawler:fractal}: in particular the set of \emph{pioneer points} (such
that $B_t$ lies on the frontier of the infinite component of the complement of
$B_{[0,t]}$),  related to the  disconnexion exponent  $\eta_1$; \emph{frontier
  points} (points of the boundary  of the infinite component of the complement
of $B_{[0,1]}$),  related to the disconnection  exponent for two  paths in the
plane.   Another exceptional  subset  of the  path  is the  set of  \emph{cone
  points} (such  that $B_{[0,t]}$ is contained  in a cone  of endpoint $B_t$),
related  to  the \emph{cone  exponents}  (studied in \cite{legall:sflour}  for
example).

We will  use the exponent  introduced in the  previous sections to  describe a
family of  exceptional sets,  indexed by  a subset $A$  of the  complex plane,
having  dimension   $2-\xi(A)$,  and   that  are  invariant   under  conformal
transformations, as follows.  Fix a Brownian path $B_{[0,1]}$, a subset $A$ of
the complex plane,  and introduce the following times  for all $t\in(0,1)$ and
$r>0$:
$$T_r(t) = \Inf\{s>t:|B_s-B_t|=r\}, \qquad S_r(t) = \Sup\{s<t:|B_s-B_t|=r\}.$$

\begin{df}
  If $0<\varepsilon<R$ and $t\in(0,1)$, let
  $$Z_t^{[\varepsilon,R]}(B)  =  \left\{  \frac {B_s-B_t}{B_{s'}-B_t}:  s  \in
    [T_\varepsilon(t), T_R(t)], s'  \in [S_R(t), S_\varepsilon(t)] \right\};$$
  and introduce $\EA^{[\varepsilon, R]}=\{B_t: Z_t^{[\varepsilon, R]} \cap A =
  \varnothing\}$. Then, letting $\varepsilon$ go to $0$:
  $$Z_t^R  =  \incr_{\varepsilon  >0}  Z_t^{[\varepsilon, R]},  \qquad  Z_t  =
  \decr_{R>0}  Z_t^R,  \qquad  \tilde  Z_t =  \decr_{R>0}  \overline{Z_t^R};$$
  define  $\EA^R$, $\EA$  and $\EEA$  accordingly. 
  
  We  shall  also use  the  notation  $\TA=\{t:B_t\in\EA\}$,  for the  set  of
  \emph{$A$-exceptional  times}, and $\TTA=\{t:B_t\in\EEA\}$,  for the  set of
  \emph{$A$-strongly     exceptional      times}.
\end{df}

Note that, since $0$ is polar  for planar Brownian motion, $Z$ is well-defined
for almost any $t$. For $A=\{1\}$,  $\EA$ is the set of local cut-points; more
generally, $t$ is in $\EA$ if, and only if, for some $\varepsilon>0$, we have
$$(B_{(t, t + \varepsilon]} - B_t) \cap  A.( B_{[t - \varepsilon ,t)} - B_t) =
\varnothing,$$
so the setup looks similar to  the definition of the exponent $\xi(A)$.  It is
easy to  see that  for all fixed  $t>0$, a.s.\  $Z_t=\CC^*$ and $\tilde  Z_t =
\CC$, so that for  $A\neq\varnothing$, $P(t\in\TA)=0$, leading to $E(\mu(\TA))
=0$ \ie\ $\mu(\TA)=0$ almost surely --- hence the term ``exceptional points''.

The set  $\EA$ of \emph{$A$-exceptional  points} is generally  not conformally
invariant.  However, it is the case for strongly exceptional points:

\begin{prp}
  \label{prp:invar}
  Let  $\Phi$ be  a conformal  map on  a neighbourhood  $\Omega$ of  $0$, with
  $\Phi(0)=0$,  and let  $B^\Omega$ be  $B$ stopped  at its  first  hitting of
  $\partial\Omega$.   By  conformal  invariance  of  planar  Brownian  motion,
  $\Phi(B^\Omega)$  is  a  Brownian  path  stopped at  its  first  hitting  of
  $\partial\Phi(\Omega)$.  Moreover, we have
  $$\EEA(\Phi(B^\Omega)) = \Phi(\EEA(B^\Omega)).$$
\end{prp}

\begin{dem}
  We  prove that  $\tilde  Z$ is  invariant.  It is  sufficient  to prove  the
  following characterization:
  $$z\in\tilde  Z_t(B) \iff  \exists  (s_n)\downarrow0,\; (s'_n)\downarrow0:\;
  \frac {B_{t+s_n}-B_t}  {B_{t-s'_n}-B_t} \to z,$$  
  as conformal maps conserve the limits  of such quotients. Such a sequence is
  easily constructed using the very definition of $\tilde Z$.
\end{dem}

Note that nothing in the preceding uses  the fact that $B$ be a Brownian path,
except for the remark about $P(t\in\TA)$. The remaining of the present section
is dedicated to deriving the Hausdorff  dimension of $\EA$ and $\EEA$. It will
be more convenient to work in the time set, so introduce
$$\TA^{[\varepsilon,R]}  = \{  t\in[0,1] :  B^1_{[t-R,t-\varepsilon]} \cap  A. 
B^2_{[t+\varepsilon,t+R]} = \varnothing\}.$$
The  scaling property  of Brownian  motion can  then be  used to  show,  as in
\cite[lemmas   3.14--3.16]{lawler:hausdorff},   that  Theorem \ref{thm:strong}
implies the following:
\begin{equation}
  \label{eq:xisurdeux}
  P(t\in\TA   {[\varepsilon,R]})  \asymp   \left(\frac   \varepsilon  R\right)
  ^{\xi(A)/2}.
\end{equation}

\subsection{Second moments}

Fix  $R>0$.  The purpose  of this  subsection is  to give  an estimate  of the
probability that  two times $t$ and  $t'$ are $A$-exceptional  times, \ie\ are
both in  $\TA^{[\varepsilon,R]}$. To get  an upper bound on  this probability,
the idea will be to  dissociate the microscopic and macroscopic scales, giving
respectively the first and second factor in the following estimate:
$$P(t,t'\in\TA^{[\varepsilon, R]}) \leqslant c \left[\frac \varepsilon R\right]
^{\xi(A)} \left[1 \wedge |t-t'|^{-\xi(A)/2}\right].$$
If $t$ and $t'$ are  two times, introduce the ``mesoscopic'' scale $d=|t'-t|$,
and separate the following three cases:
\begin{itemize}
\item  If  $d>2R$ (long-range  interaction),  the events  $E_t\eqd\{t\in\TA^{[
    \varepsilon,  R]}\}$ and $E_{t'}$  are independent,  leading to  the right
  second-order moment;
\item If $d<2R/3$  (short-range interaction), then $E_t$ and  $E_{t'}$ lead to
  three events involving disjoint subsets of the path: $t\in\TA^{[\varepsilon,
    d/2]}$,  $t'\in\TA^{[\varepsilon,  d/2]}$  and  $t\in\TA^{[3d/2,R]}$  thus
  leading to the following bound:
  $$P(E_t,E_{t'}) \leqslant C \left( \frac {2\varepsilon}{d} \right) ^{\xi(A)}
  \left(  \frac  {3d}{2R} \right)  ^{\xi(A)/2}$$  
  (in fact those events are  \emph{not} independent; however the dependence is
  only through the positions of $B$ at fixed times, so if the mesoscopic radii
  are chosen as  $(1-\varepsilon)d$ and $(3+\varepsilon)d/2$ respectively, for
  some $\varepsilon>0$,  this dependence only  contributes up to a  constant). 
  Considering $R$ as a constant we get precisely the needed estimate;
\item  Lastly, if  $2R/3<d<2R$ (medium-range  interaction), the  trivial bound
  $P(E_x,E_y) \leqslant C  (2\varepsilon/d)^{2\xi(A)}$ (obtained by forgetting
  what happens after radius $d/2$) gives the needed contribution.
\end{itemize}

So in the  case of exceptional points defined locally,  bounds
on second moments are not
difficult  to  derive (and  this  ``scale separation''  can  be  used in  various
setups).  In contrast, if the whole  path was to influence every single point,
interactions would not be that easy to classify.

\subsection{Hausdorff dimensions}

The main result of this section is the following:

\begin{thm}
  \label{thm:hausdorff}
  Let $(B_t)_{t\in[0,1]}$ be a planar Brownian path. If $A$ is any nice subset
  of the complex plane such that $\xi (A) \le 2$, then almost surely
  $$\dim_H(\EA(B)) = \dim_H(\EEA(B)) = 2-\xi(A).$$  
  In particular,  both subsets are  a.s.\ non-empty and  dense in the  path if
  $\xi(A)<2$.  If $\xi (A) >2$, $\EA(B)=\EEA(B) = \emptyset$ almost surely.
\end{thm}

\begin{dem}
  The first step in the proof is the statement of a zero-one law:

  \begin{lem}
    \label{lem:ps}
    The  dimension  of  the  set  of all  $A$-exceptional  points  (resp.\  of
    $A$-strong exceptional  points) has an almost sure  value. More precisely,
    there exist $\delta_A$ and $\tilde\delta_A$ in $[0,2]$ such that
    $$P(\dim_H(\EA)=\delta_A)=P(\dim_H(\EEA)=\tilde\delta_A)=1.$$          
    Moreover,  the following  holds with  probability  $1$ (and  the same  for
    $\EEA$ also):
    $$\forall s<t \quad \dim_H(\EA(B_{[s,t]})) = \delta_A.$$
  \end{lem}

  \begin{dem}
    The proof is the same in both cases; we perform it here for $\delta_A$.
    
    Introduce  the  following random  variables  in $[0,2]$:  $Z=\dim_H(\EA)$,
    $Z_-=\dim_H(\EA(B_{[0,1/3]}))$,     $Z_+=\dim_H(\EA(B_{[2/3,1]}))$.    The
    scaling property,  associated with the  Markov property, shows  that these
    three  variables have  the same  law;  basic properties  of the  Hausdorff
    dimension imply  that $Z\geqslant Z_-\vee  Z_+$; and locality  proves that
    $Z_-$ and $Z_+$ are independent.
    
    $0\leqslant Z_-\leqslant Z\leqslant 2$ with the same mean value: from here
    follows  that  $P(Z_-=Z)=1$.  By  the same  argument  $P(Z_+=Z)=1$,  hence
    $P(Z_-=Z_+)=1$; $Z_-$  and $Z_+$ being independent, this  is only possible
    if  they are  deterministic: thus  giving the  existence of  $\delta_A$ as
    their common almost sure value.
    
    Now if  $0\leqslant s<t\leqslant 1$  the dimension of  $\EA(B_{[s,t]})$ is
    (almost surely) $\delta_A$.  This holds  at the same time for all rational
    $s$, $t$; then  it suffices to note that  $\dim_H(\EA(B_I))$ is increasing
    in $I$ to extend the equality to all $s<t$.
  \end{dem}
  
  From this lemma  follows that as soon as $\EA$ has  positive dimension it is
  dense in the path.
  
  For convenience  we will  prove the result  in the  time set, \ie\  we shall
  compute  the dimension of  $\TA$; it  is known  that planar  Brownian motion
  doubles  Hausdorff dimensions  (\ie\  with probability  $1$,  for any  Borel
  subset      $I$      of      $[0,1]$,      $\dim_H(B_I)=2\dim_H(I)$      ---
  cf.~\cite{kaufman:metrique}),  whence  $\dim_H(\EA)=2\dim_H(\TA)$. Moreover,
  to avoid problems near $0$ and $1$  we shall suppose that $B$ is defined for
  $t\in\RR$ --- this will not change $\TA$ since the definition is local.
  
  \textbf{First step: lower  bound.} Fix $R>0$ and let  $A_n$ be the following
  set:
  $$A_n  = \{t :  B_{[t-R, t-2^{-n}]}-B_t)  \cap A(B_{[t+2^{-n},  t+R]}-B_t) =
  \varnothing\}.$$  
  For shorter notations,  let $s=\xi(A)/2$; moreover, assume from  now on that
  $s\in(0,1)$ (if $s\geqslant  1$ there is nothing to  prove, and since $A\neq
  \varnothing$ we have $s>0$ anyway).   From the previous estimates for first-
  and second-moments, we obtain
  $$E(\un_{A_n}(x)) \asymp 2^{-sn} \qquad E(\un_{A_n}(x) \; \un_{A_n}(y))
  \leqslant c 2^{-sn} \left[ 1 \wedge \frac {2^{-sn}} {|y-x|^s} \right].$$  
  Introduce the (random) measure $\mu_n$ having density $2^{sn}\un_{A_n}$ with
  respect to  the Lebesgue  measure. It  is not hard  to derive  the following
  estimates:
  \begin{align}
    \label{al:m1} E(\|\mu_n\|) &= \int_{[0,1]} 2^{sn} E(\un_{A_n}(x)) \;\dd x
    \asymp 1, \\
    \label{al:m2} E(\|\mu_n\|^2) &= \int\!\!\!\!\int_{[0,1]^2} 2^{2sn}
    E(\un_{A_n}(x) \; \un_{A_n}(y)) \;\dd x\,\dd y \nonumber \\
    &\leqslant c2^{sn} \left[\int_0^1\dd x\int_x^{x+2^{-n}} \dd y \;+\; \int
    _0 ^{1-2^{-n}} \dd x \int_{x+2^{-n}}^1 \frac{2^{-sn}\dd y} {(y-x)^s}
    \right] \nonumber \\
    &\leqslant c2^{(s-1)n} \;+\; c \int _0 ^{1-2^{-n}} \left( \frac
    {(1-x)^{1-s}} {1-s} - \frac {2^{(s-1)n}} {1-s}\right) \dd x \nonumber \\
    &\leqslant c + c2^{(s-1)n} + c2^{(s-2)n} \leqslant c.
  \end{align}
  Hence, $\|\mu_n\|$  has finite expectation and  finite variance, independent
  of $n$: there exists  $\varepsilon>0$ satisfying $P(\|\mu_n\| > \varepsilon)
  > \varepsilon$  for all positive  $n$.  Consequently,  it is  possible, with
  positive probability, to extract  a subsequence $(\mu_{n_k})$ such that, for
  all $k$,  $\|\mu_{n_k}\| \geqslant \varepsilon$. By  a compactness argument,
  another extraction  leads to  a converging subsequence,  the limit  $\mu$ of
  which satisfies  $\|\mu\| \geqslant \varepsilon$. $\mu$ is  supported on the
  intersection of the $A_n$,  this intersection is non-empty: hence $P(\bigcap
  A_n\neq\varnothing)>0$.
  
  Introduce then the notion of $r$-energy  of a measure: if $\nu$ is some mass
  measure supported on a metric space $X$, let
  $$\mathcal E_r(\nu) \eqd \int\!\!\!\!\int_{X^2} \frac {\dd\nu(x)\,\dd\nu(y)}
  {\dd(x,y)^r}.$$  
  It is known  that if $X$ supports a mass measure  of finite $r$-energy, then
  its Hausdorff dimension is not less than $r$ (cf.~\cite{kahane:random}). Let
  then   $r\in(0,1-s)$:   a    calculation   analogous   to   the   derivation
  of (\ref{al:m2})  leads to
  \begin{equation}
    \label{eq:energy}
    E(\mathcal  E_r(\mu_n)) \leqslant  c + c2^{(r+s-1)n} + c2^{(r+s-2)n}
    \leqslant c.
  \end{equation}
  Performing another  subsequence extraction, it  is possible to  obtain $\mu$
  supported on $\bigcap A_n$ and having finite $r$-energy: hence
  $$\forall r<1-s\quad P(\dim_H(\bigcap A_n)  \geqslant r) >0.$$  
  By  definition $\TA$  is the  increasing  union, for  $R$ going  to $0$,  of
  $\bigcap_n A_n(R)$:  hence for all $r<1-s$ we  have $P(\dim_H(\TA) \geqslant
  r)>0$.   Combining this  and the  zero-one result  (Lemma~\ref{lem:ps}) then
  proves that almost surely $\dim_H(\TA)\geqslant 1-s$.
  
  \textbf{Second step: upper bound.} This  step is usually the easier one, but
  in  the  present  case a  complication  arises  due  to  the fact  that  the
  ``non-intersection'' event we  consider at $B_t$ depends on  the position of
  $B_t$  ---   which  is   not  the   case  for  instance   in  the   case  of
  cut-points \cite{lawler:fractal}.  This  explains   why  we  need  one  more
  argument, namely the continuity of $\xi:A\mapsto\xi(A)$.
  
  Fix   a   nice   set    $A$,   $\varepsilon>0$,   $R>0$   and   a   sequence
  $(\lambda_n)_{n\geqslant0}$ of  positive numbers, tending slowly  to $0$ (in
  the following sense:  for all positive $\eta$, $2^{-\eta  n} = o(\lambda_n)$
  --- for instance,  take $\lambda_n=1/n$).  Now  suppose some time $t$  is in
  $A_n$. With positive probability, the following happens:
  $$\left\{\begin{array}{l}
      B_{[t-\lambda_n 2^{-n},t+\lambda_n 2^{-n}]} \subset \mathcal
             B(B_t,\lambda_n^{1/2} 2^{-n/2})\\ 
      |B_{t-2^{-n}}-B_t| \geqslant 2^{-n/2}\\
      |B_{t+2^{-n}}-B_t| \geqslant 2^{-n/2}\\
      (B_{[t-R,t-2^{-n}]} \cup B_{[t+2^{-n},t+R]}) \cap \mathcal
      B(B_t,(1-\varepsilon)2^{-n/2}) = \varnothing
    \end{array}\right.$$  
  (the first three conditions are a consequence of scaling, and the fourth one
  is the start-separation lemma, more  precisely the weakened version of it as
  stated in equation (\ref{eq:nowedge})). Introduce  $A^{\eta_n} = \{az : a\in
  A, z\in\mathcal B(1,\eta_n)$: we have
  $$P(B_{[t-R,  t-2^{-n}]}-B_t) \cap  A^{\eta_n}  (B_{[t+2^{-n}, t+R]}-B_t)  =
  \varnothing \;|\; t\in A_n)$$
  \nopagebreak\vspace{-2em}
  \begin{equation}
    \label{eq:enlarge}
    \asymp \frac {2^{-n\xi(A^{\eta_n})/2}} {2^{-n\xi(A)/2}} =
    2^{-n[\xi(A^{\eta_n})-\xi(A)]/2}.
  \end{equation}
  It is  easy to  see that  under the previous  conditions, if  $t\in \mathcal
  T_{A^{\eta_n}}$,  then every $t'\in[t-\lambda_n  2^{-n},t+\lambda_n 2^{-n}]$
  is in $A_n$, as soon  as $\eta_n > 18\lambda_n/(1-\varepsilon)$. From now on
  we  shall assume  that  this  holds, and  that  $\eta_n\to0$. Putting  these
  estimates  together, we  obtain the  following  (where $l$  is the  Lebesgue
  measure on $\RR$): for all interval $I$,
  \begin{equation}
    \label{eq:exists}
    P(l(A_n\cap I)>\lambda_n 2^{-n} | A_n\cap I \neq\varnothing) \geqslant
    c.2^{-n[\xi(A^{\eta_n})-\xi(A)]/2}.
  \end{equation}

  The Markov inequality then states that
  $$P(l(A_n\cap I)>\lambda_n 2^{-n}) \leqslant \frac {E(l(A_n\cap I))}
  {\lambda_n 2^{-n}},$$   
  and   $E(l(A_n\cap    I))   \asymp   2^{-n\xi(A)/2}    l(I)$.    From   this
  and (\ref{eq:exists}) follows that
  \begin{equation}
    \label{eq:exists2}
    P(A_n\cap I \neq\varnothing) \leqslant C\, \frac {2^{-n\xi(A)/2} l(I)}
    {\lambda_n 2^{-n}} \,\frac 1 {2^{-n[\xi(A^{\eta_n})-\xi(A)]/2}}.
  \end{equation}
  By continuity  of $\xi$, for  large $n$ we have  $|\xi(A^{\eta_n})-\xi(A)| <
  \varepsilon$; by the hypothesis on  $\lambda_n$, still for large $n$ we have
  $\lambda_n \geqslant 2^{-\varepsilon/2}$. Hence for large $n$:
  \begin{equation}
    \label{eq:majo}
    P(A_n\cap I \neq\varnothing) \leqslant C\, 2^{\varepsilon n}\,
    2^{-n\xi(A)/2} \, \frac{l(I)}{2^{-n}}.
  \end{equation}
  Cover the  interval $[0,1]$ with the  $I_k^n=[k2^{-n},(k+1)2^{-n}]$, and let
  $X_n$ be the number of such intervals intersecting $A_n$. Then
  \begin{align*}
    E(X_n) &= \sum_k P(I_k^n\cap\TA \neq \varnothing) \leqslant 2^n\; C\,
    2^{\varepsilon n}\, 2^{-n\xi(A)/2} \, \frac{l(I\z^n)}{2^{-n}} \leqslant C
    \, 2^{\varepsilon n} \, 2^{n[1-\xi(A)/2]}.
  \end{align*}  
  By another application of the Markov inequality,
  $$P(X_n>2^{n[1-\xi(A)+2\varepsilon]})   \leqslant  C\,  2^{-\varepsilon n}.$$
  Hence by  the Borel-Cantelli theorem,  for sufficiently large $n$,  $A_n$ is
  covered  by  at  most  $2^{n[1-\xi(A)+2\varepsilon]}$  intervals  of  length
  $2^{-n}$  ---   and  this   implies  that  $\dim_H(\bigcap   A_n)  \leqslant
  1-\xi(A)/2+2\varepsilon$. Letting  $\varepsilon$ tend  to $0$ then  leads to
  (a.s.)  $\dim_H(\bigcap  A_n) \leqslant 1-\xi(A)/2$.   This is true  for all
  $R>0$, hence remains true in the limit $R\to0$: together with the first step
  of  the proof  this  gives (a.s.)   $\dim(\TA)=1-\xi(A)/2$ hence  $\dim(\EA)
  =2-\xi(A)$.
  
  Then, $\EEA$ is contained in  $\EA$ and contains every $\mathcal E_{A^\eta}$
  for  positive $\eta$  (with  the  previous notations):  another  use of  the
  continuity of $\xi$ then gives $\dim_H(\EEA)=\dim_H(\EA)=2-\xi(A)$.
\end{dem}

As a consequence, we get a second result:

\begin{thm}
  If  $A$  is  any  nice  subset  of  the  complex  plane,  then  the  set  of
  \emph{globally $A$-exceptional points}, \ie\ points $B_t$ satisfying
  $$(B_{[0,t)} - B_t) \cap A.(B_{(t,1]}  - B_t) = \varnothing,$$  
  has  Hausdorff  dimension $2-\xi(A)$  ---  and  in  particular it  is  a.s.\ 
  non-empty for $\xi(A)<2$, and a.s.\ empty for $\xi(A)>2$.
\end{thm}

\begin{dem}
  Again, extend $B$ to $(B_t)_{t\in\RR}$ defined on the entire real line.  The
  set  $\TA^1$  of  $A$-exceptional  times  up  to the  scale  $R=1$  (as  was
  introduced previously) in $[0,1]$ is exactly the set of globally exceptional
  points. Therefore,  the previous proof  can be applied directly.   The upper
  bound  is  immediate: since  every  globally  exceptional  point is  locally
  exceptional we  have $\dim_H(\TA^1)\leqslant\dim_H(\TA)\leqslant 1-\xi(A)/2$
  a.s.
  
  The  lower bound  requires  a little  more work,  indeed  we do  not have  a
  zero-one  law for the  dimension of  $\TA^1$. It  can be  seen that  in fact
  Equation~(\ref{eq:energy}) can be refined, the proof being exactly the same,
  into the following (with the same notations as previously):
  $$\exists C>0 \quad \forall  r\in(0,1-s) \quad \forall n>0 \qquad E(\mathcal
  E_r(\mu_n)) \leqslant \frac C {1-(r+s)},$$  
  where $C$ may only depend on $A$.  Hence, with the same constant and for all
  $\lambda>1$:
  $$P\left(\mathcal  E_r(\mu_n) \leqslant  \frac {\lambda  C} {1-(r+s)}\right)
  \geqslant 1-\frac1\lambda.$$  
  one   can   then   perform   the  subsequence   extraction   (cf. proof   of
  Theorem \ref{thm:hausdorff}) in a way which ensures that, for all $r$,
  \begin{equation}
    \label{eq:bilan}
    P\left(\|\mu\|>0 \;\mbox{and}\; E_r(\mu) \leqslant  \frac {\lambda C}
    {1-(r+s)} \right) \geqslant c,
  \end{equation}
  with $c>0$ and $\lambda>1$ independent of $r$. Moreover, $\mathcal E_r(\mu)$
  being a non-decreasing function of $r$ (since the set $[0,1]$ is of diameter
  $1$), we  finally obtain,  with positive probability,  a mass  measure $\mu$
  supported on $\TA$ satisfying
  $$\forall  r<1-s  \qquad  \mathcal  E_r(\mu)  \leqslant  \frac  {\lambda  C}
  {1-(r+s)} <\infty.$$    
  Hence, with positive probability,  $\dim_H(\TA) \geqslant 1-s = 1-\xi(A)/2$,
  and combining this to the previous paragraph leads to
  $$P\left(\dim_H(\TA) = 1-\frac{\xi(A)}2 \right) > 0.$$
  It   is   then   possible   to   conclude   using   the   same   method   as
  in \cite[pp. 8--9]{lawler:hausdorff}.
\end{dem}

\subsection{Remark about critical cases}

In cases  where $\xi(A)=2$, the previous  theorem is not  sufficient to decide
whether $A$-exceptional points exist. We  shall see in the next paragraph that
$\xi((-\infty,0))  = \xi((0,\infty))=2$.  In  fact these  two  cases are  very
different:

\begin{prp}
  Almost  surely, $\EA$  is  empty for  $A=((0,\infty))$  and non-empty  (with
  Hausdorff dimension $0$ though) for $A=((-\infty,0))$.
\end{prp}

\begin{dem}
  The second point is easier: if $t$ is such that $\Re(B_t)$ is maximal in the
  path, then  $B_{[0,1]}$ lies inside  a half-plane whose border  goes through
  $B_t$.  Since a.s.\  $B_t$ is  the only  point having  this real  part, this
  proves  that $(B_s-B_t)/(B_{s'}-B_t)$  is never  in $(-\infty,0)$,  which is
  precisely what we wanted.
  
  The first point is more problematic.  The method used to derive the value of
  $\xi$ for a  wedge with end-point at the  origin (cf. next paragraph) allows
  to prove the following: Let $\alpha$  and $\beta$ be in $(0,2\pi)$, then the
  probability that, given independent paths  $B^1$ and $B^2$ starting from the
  unit  circle, there exist  two wedges  of angles  $\alpha$ and  $\beta$, and
  containing respectively $B^1$ and $B^2$ up to radius $R$, decreases as
  $$p_R(\alpha,\beta) \approx R^{\pi/\alpha + \pi/{\beta}}.$$  
  Hence, as  soon as  $\pi/\alpha+\pi/{\beta}$ is greater  than $2$,  there is
  a.s.\ no point  $B_t$ on the path  such that $B_{[0,t]}$ lies in  a wedge of
  angle $\alpha$ and $B_{[t,1]}$ lies in a wedge of angle $\beta$ (there is no
  ``asymmetric two-sided cone point'' of those angles on the path).
  
  For all $\alpha\in(0,\pi)$,  introduce $\alpha_1=2\pi-\alpha$ and $\alpha_2$
  as the  biggest angle in $(0,2\pi]$  satisfying $\pi/\alpha+\pi/\alpha_2>2$. 
  Note that $\alpha_2>\alpha_1$: denote then
  $$\beta(\alpha)=\frac{\alpha_1+\alpha_2}{2}.$$  
  Note  that $\pi/\alpha+\pi/\beta(\alpha)>2$  and $\beta(\alpha)+\alpha>2\pi$
  for all $\alpha\in(0,\pi)$.  From this follows that, almost  surely, for all
  $\alpha\in(0,\pi)\cap\QQ$,  there is  no assymetric  cone point  with angles
  $\alpha$ and $\beta(\alpha)$.
  
  Let now $A=(0,\infty)$ and suppose there is a point $B_t$ in $\EA$. That is,
  there  exist two  half-lines  starting from  $B_t$  whose reunion  separates
  $B_{[0,t]}$ from $B_{[t,1]}$. Then we are in one of two cases:
  \begin{itemize}
  \item Either these half-lines form a straight line, \ie\ there is a straight
    line cutting the path. This cannot  happen, as recently proved by Bass and
    Burdzy \cite{bass:cutting} --- and the proof is very difficult.
  \item  Or  there  are  disjoint  wedges  of  angles  $\alpha\in(0,\pi)$  and
    $2\pi-\alpha$, each  containing one part  of the path. Then,  there exists
    $\alpha\z\in\QQ$   such  that   $\alpha\z>\alpha$   and  $\beta(\alpha\z)>
    2\pi-\alpha$, and $B_t$ is an asymmetric cone point with angles $\alpha\z$
    and $\beta(\alpha\z)$. We just saw that such a point cannot exist.
  \end{itemize}
  Hence $\EA=\varnothing$.
\end{dem}


\section{Bounds and conjectures on the exponent function}
\label{sec:numerical}

\subsection{Known exact values of $\xi$}

\begin{prp}
  \begin{enumerate}
  \item $\xi(\{1\})=5/4$, hence for all $z\neq0$ and $n>0$:
    $$\xi \left( \{ze^{2ik\pi/n},k=1,\ldots,n\} \right) =5n/4;$$
  \item Letting $W_\alpha$ be a wedge of angle $0\leqslant \alpha <2\pi$:
    $$\xi(W_\alpha) = \frac {4\pi} {2\pi-\alpha};$$
    in particular $\xi((0,\infty))=\xi((-\infty,0))=2$;
  \end{enumerate}
\end{prp}

\begin{dem}
  (i): The value of $\xi (\{1 \})  = 5/4$ has recently been derived by Lawler,
  Schramm and  Werner \cite{werner:value2},  and the proof  is far  beyond the
  scope  of this  paper. The  result  for all  $n$ is  then a  straightforward
  consequence of Proposition \ref{prp:props}, point (iv).
  
  (ii): Suppose  $A=W_\alpha$ is  centered around the  positive axis,  so that
  $A=\{re^{i\theta},  r>0, |\theta| <  \alpha/2\}$; introduce  the symmetrical
  wedges  $W_\beta'=\{re^{i\theta}, r>0, |\theta-\pi|  < \beta/2\}$.  If $B^1$
  stays in  $W_{\pi-\alpha/2}$ and $B^2$ remains  in $W'_{\pi-\alpha/2}$, then
  $B^1\cap AB^2=\varnothing$: The  probability of staying in a  wedge of angle
  $\beta$ until  radius $R$ being  strongly approximated by  $R^{-\pi/ \beta}$
  (the exponent is obtained through  the gambler's ruin estimate combined with
  the  analyticity of the  exponential function;  the strong  approximation is
  true but  in fact not  needed here, cf.  \cite{evans:cone}), we get  a lower
  bound:
  $$p_R(W_\alpha) \geqslant c \left(R^{-\pi/(\pi-\alpha/2)}\right)^2,$$
  hence $\xi(W_\alpha) \leqslant 4\pi/(2\pi-\alpha)$.
  
  Now  remark that  the condition  $B^1\cap AB^2=\varnothing$  means  that the
  complement of  the paths  contains an ``hourglass'',  \ie\ the union  of two
  disjoint wedges  of angle $\alpha/2$.  So introduce $\eta>0$ and  a (finite)
  family  $(S_i)_{1\leqslant   i\leqslant  N}$  of   hourglasses  with  angles
  $\alpha/2-\eta$, such that any  hourglass with angle $\alpha/2$ contains one
  of the  $S_i$. If $q_R(i)$ is  the probability that the  paths are separated
  from  each other  by  $S_i$,  then $p_R(W_\alpha)  \leqslant  \sum q_R(i)$.  
  Noticing  that if  $\beta_i$ and  $\beta'_i$ are  the angles  of  the wedges
  forming the complement of $S_i$,  we obtain as previously $q_R(i) \asymp R^{
    -\pi/\beta_i  -\pi/\beta'_i}$, and  optimizing this  under  the constraint
  $\beta_i+\beta'_i=2\pi-(\alpha-2\eta)$ ---  where the greatest  value is for
  $\beta=\beta'$ --- we finally get the following estimate:
  $$p_R(W_\alpha)  \leqslant  CN\,R^{-2\pi/(\pi+\eta-\alpha/2)}.$$  
  From  this follows  that  $\xi(W_\alpha)\geqslant 4\pi/(2\pi+2\eta-\alpha)$,
  and letting  $\eta$ go  to $0$ then  gives the  conclusion --- at  least for
  $\alpha>0$. But in fact the  same method still applies for $\alpha \geqslant
  0$: simply  inflate the complement  of the hourglass instead  of introducing
  angle $\alpha/2-\eta$, the fact that the wedges to consider may overlap does
  not change anything to the proof.
\end{dem}

\emph{Remark:} If we denote $A^\alpha=\{ze^{i\theta},z\in A,|\theta| \leqslant
\alpha/2\}$ (that  is, $A$  ``thickened'' by an  angle $\alpha$), then  it can
easily be proven that
\begin{equation}
  \label{eq:bound}
  \xi(A^\alpha) = \frac {h_A(\alpha)} {2\pi-\alpha},
\end{equation}
where  $h_A$ is  continuous  (until the  angle  $\alpha\z\leqslant 2\pi$  when
$\xi(A^\alpha)$   tends   to    infinity),   non-decreasing,   and   satisfies
$h_A(0)=2\pi\xi(A)$; in the wedge case, $h$ is constant.

\subsection{An upper bound for the exponent}
\label{sub:bound}

From  continuity of  $\xi$ and  the  exact value  $\xi(\{1\})=5/4<2$, one  can
deduce that there  are ``pivoting points'' of any  sufficiently small angle on
the  Brownian path (that  is, points  around which  one half  of the  path can
rotate of a small angle without intersecting the other half --- the associated
$A$ being  $\mathcal C_\alpha = \{e^{i\theta}, \theta  \in [0,\alpha]\}$). The
following proposition gives a (bad  but) quantitative bound for such values of
$\alpha$ --- without usage of the exact value for $\alpha=0$:

\begin{prp}
  \label{prp:bound}
  For all positive $\alpha$, we have the following upper bound:
  $$\xi(\mathcal C_\alpha) \leqslant \frac{4\pi}{2\pi-\alpha} \left[1 - \frac
    {(\log2)^2} {4\pi^2}\right].$$    
\end{prp}

\begin{dem}
  The proof is adapted  from \cite{werner:disconnection}, where an upper bound
  for the  classical disconnection exponent for one  path, \ie $\xi(1,0)$, was
  obtained. The method  is the following: First, estimate  the extremal length
  of  a strip  bounded by  Lipschitz functions;  then describe  a sufficiently
  large subset of  $E_R$, using such strips, and use  the previous estimate to
  derive a bound for $P(E_R)$.

  \begin{inlem}
    Let  $f$ be  a  continuous, $M$-Lipschitz  function  on $\RR$,  satisfying
    $f(x)+f(-x)=2f(0)$ for all $x$, and  let $\beta>0$. Introduce the strip of
    width $\beta$ and length $2r$ around $f$ as
    $$\mathcal B_f^\beta(r) = \left\{ x+iy \;:\; |x| < r,\; |y-f(x)| <
      \frac\beta2 \right\};$$        
    let  $W$  be  a planar  Brownian  path  starting  at $if(0)$,  and  denote
    \smash{$A_f^\beta (r)$}  the event that  the point $x+iy$ where  $W$ first
    reaches $\partial \mathcal  B_f^\beta(r)$ satisfies $|x|=r$ (\ie $W$ exits
    $\mathcal B$ by one of the vertical parts of its boundary). Then
    $$P(A_f^\beta(r)) \geqslant \frac1\pi \exp\left[ -\frac{\pi r}{\beta}
      (1+M^2) \right].$$
  \end{inlem}
  
  \begin{dem}
    This is an easy consequence of  the following estimate, which can be found
    in \cite{ahlfors:confinv}      and      is      a      consequence      of
    Proposition \ref{prp:extbound2}: If  $L$ is the  extremal distance between
    both vertical parts of $\partial\mathcal B$ in $\mathcal B$, then
    $$L\leqslant \frac {2r}\beta (1+M²);$$    
    using this together with the classical estimate for Brownian motion in a
    strip provides the right estimate.
  \end{dem}
  
  For  the rest  of this  proof, we  shall consider  paths in  the logarithmic
  space, denoted by  the letter $W$; the actual path $B$  is obtained from $W$
  by applying the exponential map  --- conformal invariance of Brownian motion
  then proves that  $B$ is a Brownian path.  Let $f$ be a function  such as in
  the lemma:  it is clear that  if $W^1$ remains in  $\mathcal B_f^\pi(r)$ and
  $W²$  stays  in $\mathcal  B_{f+\pi}^\pi(r)$,  then  $B¹$  and $B²$  do  not
  intersect up to the first time they reach radius $e^r$ or $e^{-r}$. Together
  with  the  fact  that  $P(A_f^\pi(r))=P(A_{f+\pi}^\pi(r))$,  this  leads  to
  $P(E_R(\{1\})) \geqslant (P(A_f^\pi(\log R))/2)²$, hence using the lemma:
  \begin{equation}
    \label{eq:one_f}
    P(E_R(\{1\})) \geqslant c R^{-2(1+M²)}.
  \end{equation}
  Doing the same  with strips of width $\beta=\pi-\alpha/2$  (for which it can
  be  seen that  $B¹$ and  $B²$ can  rotate around  $0$ by  an angle  at least
  $\alpha/2$ in each direction) leads to
  \begin{equation}
    \label{eq:one_alpha}
    P(E_R(\mathcal C_\alpha)) \geqslant c \exp\left[
    -\frac{4\pi}{2\pi-\alpha}(1+M^2) \log R \right],
  \end{equation}
  hence, letting $f=0$, a first bound on the exponent:
  $$\xi(\mathcal C_\alpha) \leqslant \frac{4\pi}{2\pi-\alpha}$$    
  (this is also  a direct consequence of $\mathcal  C_\alpha \subset W_\alpha$
  and the  exact value of $\xi(W_\alpha)$,  which happens to  be precisely the
  upper bound we  just obtained). Note that the bound is  never less than $2$,
  hence we proved nothing useful yet.
  
  We now  want to consider  families of strips. Keep  $\beta=\pi-\alpha/2$ and
  fix  $\gamma>0$;  let  $U_N=\{±1\}^N$  and  for  $u\in  U_N$  let  $f_u$  be
  constructed as follows:
  \begin{itemize}
  \item  $f_u(0)=0$, and  for $1\leqslant  n\leqslant N$,  $\displaystyle f_u(
    n\gamma )=\frac\beta2 \sum_{k=1}^{n} u_k$;
  \item   $f$   is   affine   on   each   $[n\gamma,(n+1)\gamma]$,   satisfies
    $f_u(x)=f_u(N\gamma)$  for all $x>N\gamma$  and $f_u(-x)=-f_u(x)$  for all
    $x$.
  \end{itemize}
  Then for $u\neq u'$ the intersection of \smash{$\mathcal B_{f_u}^\beta$} and
  \smash{$\mathcal B_{f_{u'}}^\beta$} is  not connected, hence \smash{$A_{f_u}
    ^\beta$} and \smash{$A_{f_{u'}}^\beta$} are disjoint. This leads to
  $$P(E_R (\mathcal C_\alpha)) \geqslant c \sum_{u\in U_N} \exp\left[
    -\frac{2\pi}{\beta} (1+(\beta/2\gamma)^2) \log R \right]$$    
  for all  $N$, where $R=e^{N\gamma}$. Then using  $P(E_R (\mathcal C_\alpha))
  \asymp R^{-\xi(\mathcal C_\alpha)}$, noticing that  all the terms of the sum
  are equal (there are $2^N$ of them) and applying a logarithm:
  \begin{equation}
    \label{eq:many_beta_N}
    \xi(\mathcal C_\alpha)N\gamma \leqslant \frac {2\pi}
    {\beta} (1+(\beta/2\gamma)^2) N \gamma - N \log2 -\log c.
  \end{equation}
  Divide by $N\gamma$ and let $N$ go to infinity to obtain
  \begin{equation}
    \label{eq:final}
    \xi(\mathcal C_\alpha) \leqslant \frac{\pi\beta}{2}
    \left(\frac1\gamma\right)^2 - \log2 \left(\frac1\gamma\right) + \frac
    {2\pi} {\beta}.
  \end{equation}
  This  is true  for all  $\gamma>0$;  the optimal  value is  $\gamma=\pi\beta
  /\log2$, leading to
  $$\xi(\mathcal C_\alpha) \leqslant \frac {4\pi} {2\pi-\alpha} \left[ 1 -
    \frac{(\log2)^2} {4\pi^2} \right],$$
  which is precisely what we wanted.
\end{dem}

\emph{Remark:} The same proof gives a  bound on $\xi(A)$ if $A$ is included in
a small ball centered  at $1$, as a function of the  radius. But since it does
not make use of the value  of $\xi(\{1\})$, no modulus of continuity for $\xi$
can be  obtained this  way. Cf. however equation (\ref{eq:bound})  for another
bound, which does provide such a modulus but is not quantitative.

\bigskip

As a consequence of this bound, we obtain the following

\begin{thm}
  For all  $\alpha<\log^22/2\pi$, the  following holds: With  probability $1$,
  the set of local pivoting points of angle $\alpha$ on a planar Brownian path
  is non-empty and has a positive Hausdorff dimension.
\end{thm}

\begin{figure}[h!]
  \begin{center}
    \leavevmode
    \includegraphics[scale=1]{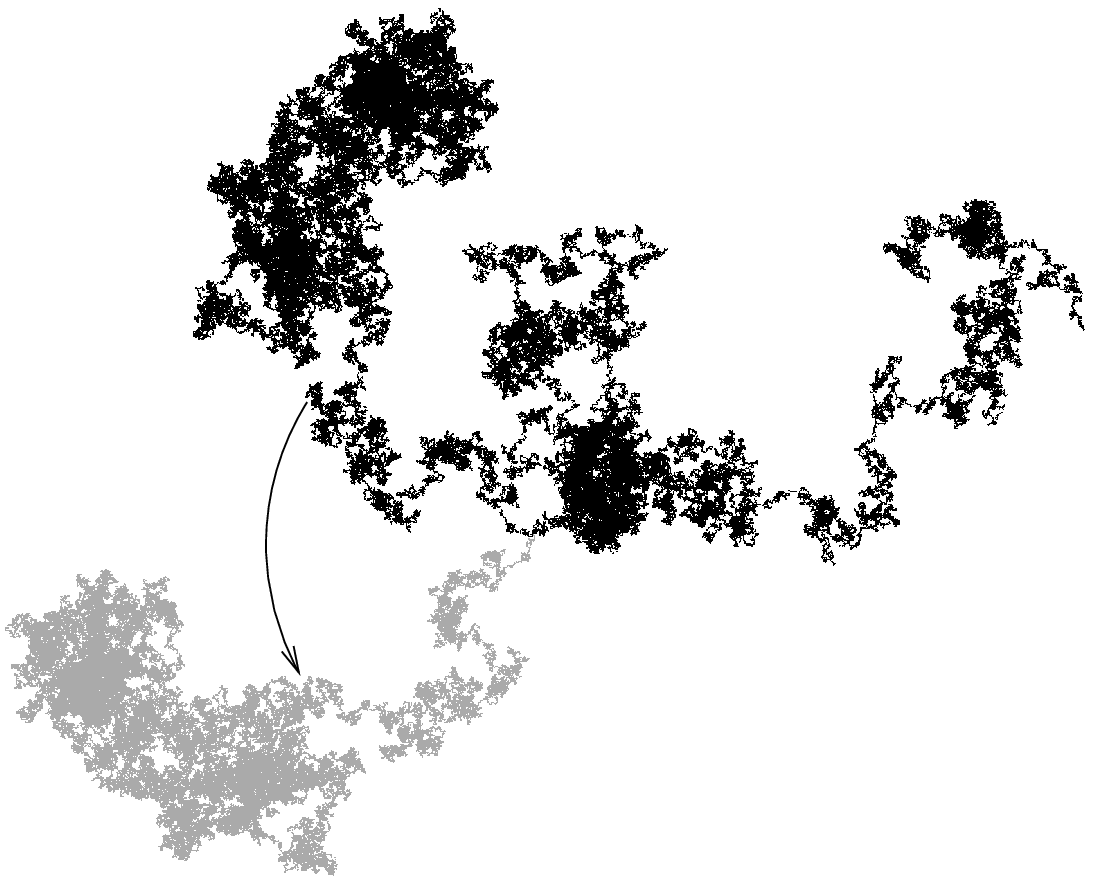}
    \caption{A pivoting point of angle $\pi/2$}
    (in grey is the image of one half of the path by a rotation of angle
    $+\pi/2$)
    \label{fig:piv}
  \end{center}
\end{figure}

\emph{Remark:} The bound given in the theorem ($\log^22/2\pi \simeq 0.076$) is
certainly not the best one; simulations suggest that there are pivoting points
of  any angle  less than  $3\pi/4  \simeq 2.356$  --- cf. next  subsection for
details and figure \ref{fig:piv} for a picture of a pivot of angle $\pi/2$. In
particular, the maximal angle is  conjectured to be greater than $2\pi/3$, and
this seems  to indicate  that a discrete  analogue of (local)  pivoting points
will appear  on the exploration process  of a critical  percolation cluster on
the triangular lattice \cite{schramm:UST,smirnov:perco}.

\subsection{Conjectured and experimental values}
\label{subsec:conj}

Some  exact values  of $\xi(A)$  are known,  cf. subsection \ref{sub:propexp}. 
However,  heuristic arguments  seem to  indicate that  the formula  giving the
exponent  for wedges is  close to  apply in  other cases  such as  notably the
``weak pivot'' exponent, namely:
$$\xi(\{1,e^{i\theta}\})    \simeq   \frac{5\pi/2}{2\pi-\theta}$$
for  all  $\theta\in[0,\pi]$ ---  corresponding  to  a  continuous version  of
Proposition \ref{prp:props}, point (iv). This  is confirmed by simulations, at
least for $\theta=\pi/2$  and $\theta=\arctg (3/4)$ (cf. table \ref{tab:exp}),
based on the following

\begin{hogtheorem}{Conjecture}
  Let $A$ be  a bounded, non-empty subset of  $\ZZ^2\setminus\{0\}$; let $B^1$
  and $B^2$ be  independent Brownian paths starting respectively  from $0$ and
  $1$,  and $S^1$  and $S^2$  be  independent standard  random walks  starting
  respectively from $0$ and $(a,0)$ with  $a$ sufficiently large (so as not to
  make the probability in the formula equal to $0$).  Then,
  $$P(B^1_{[0,T]} \cap AB^2_{[0,T]} = \varnothing) \asymp P(S^1_{[0,T]} \cap
  AS^2_{[0,T]} = \varnothing) \asymp T^{-\xi(A)/2}.$$
\end{hogtheorem}

\begin{dem}
  There  is no  known  direct proof  of  the existence  of a  non-intersection
  exponent for random  walks, the only way to obtain  the desired behaviour is
  coupling  with  Brownian  motion  --- cf. \cite{lawler:walk}.   The  present
  generalization can certainly be obtained in a similar way, note however that
  walks appear  that are not  standard simple random  walks but take  steps in
  $\{a,ia,-a,-ia\}$ for some $a \in  \CC \cap \ZZ^2$; exponents for such walks
  are the same as for SRW's (cf. \cite{cranston:rw}), but strong approximation
  in not yet proved.
\end{dem}

The most severe restriction is  the assumption that $A \subset \ZZ^2 \setminus
\{0\}$,  in  particular  simulations  cannot  (yet) be  performed  if  $A$  is
connected,  except for  very special  cases such  as wedges  (where  the exact
exponent  is   known).   However  homogeneity  can  sometimes   be  used  when
$A\in\QQ^2$ (as  for $A=\{5,4+3i\}$ which  has the same exponent  as $\{1,e^{i
  \theta}\}$ for $\theta = \arctg (3/4)$).

\begin{table}[ht]
  \begin{center}
    \begin{tabular}{|c|c|c|c|c|}
      \hline
                      & conjectured & number        & computed    & relative   \\
      $A$             & exponent    & of samples    & exponent    & error      \\
      \hline
      $\{±1\}$        & $\sim2.5$   & $2.6\,10^9$   & $2.501293$  & $+0.05 \%$ \\
      \hline
      $\{1,i\}$       & $\sim5/3$   & $3.0\,10^8$   & $1.662239$  & $-0.27 \%$ \\
                      &             &               & $1.668242*$ & $+0.09 \%$ \\
      \hline
      $\{5,4+3i\}$    & $\sim1.392679$ & $1.2\,10^6$   & $1.382311$  & $-0.74 \%$ \\
                      &             &               & $1.394610*$ & $+0.14 \%$ \\
      \hline
      $\{5,4+3i,5i\}$ & $\sim5/3$   & $1.6\,10^7$   & $1.662964$  & $-0.22 \%$ \\
                      &             &               & $1.665650*$ & $-0.06 \%$ \\
      \hline
    \end{tabular}
    \caption{\setlength{\leftskip}{2cm}\setlength{\rightskip}{2cm}%
      Some  simulated  values  of  $\xi$} \small  ($100\,000$-step  walks  ---
    exponents  marked with  a star  \break  are obtained  after a  non-rigorous
    correction)
    \label{tab:exp}
  \end{center}
\end{table}


\section{Appendix}
\label{sec:ingredients}

\subsection{Sub-additivity}

The  following  proposition   is  well  known  and  included   here  only  for
completeness (note  however that  the bounds are  not asymptotic and  that the
constants are exactly known, which is needed to derive continuity of $\xi$).

\begin{prp}[Subadditivity]
  \label{prp:subadd}
  Let $f:[1,\infty)\to(0,\infty)$ be some function such that:
  \begin{itemize}
  \item $f$ is bounded and bounded away from $0$ on any $[0,l]$, $l>0$;
  \item There exist $\varepsilon$, $A$,  $c$ and $C$ in $(0,\infty)$ such that
    for    all    $t\geqslant1$,    $ct^{-A}    \leqslant    f(t)    \leqslant
    Ct^{-\varepsilon}$;
  \item There exist $0\leqslant c_-\leqslant c_+\leqslant\infty$, at least one
    of which finite and positive, such that
    $$\forall t,t'\in[1,\infty) \quad  c_-f(t)f(t') \leqslant f(tt') \leqslant
    c_+f(t)f(t').$$
  \end{itemize}

  Then, there is a $\xi>0$ such that $f(t)\approx t^{-\xi}$. Moreover, for all
  $t\geqslant 1$,
  $$c_+^{-1}t^{-\xi} \leqslant f(t) \leqslant c_-^{-1}t^{-\xi}.$$    
  In particular,  if both $c_-$  and $c_+$ are  in $(0,\infty)$ we  get strong
  approximation: $f(t)\asymp t^{-\xi}$.
\end{prp}

\subsection{Extremal distance}

Many of  the known estimates for  exponents (apart from cases  where the exact
value in known ---  such as the exponent of a cone  here, and the intersection
exponents   in   the   half-plane   in \cite{werner:value})  come   from   the
corresponding  estimates for  Brownian  paths in  rectangles, using  conformal
invariance. The  introduction of extremal  distance generalizes the  notion of
aspect ratio  of a rectangle  and hence provides  a natural parameter  in this
process.

\begin{hogtheorem}{Theorem and Definition :}
  Let  $\Omega$ be an  open, bounded,  simply connected  subset of  $\CC$, the
  frontier of  which (oriented in  the usual direct  sense) is a  Jordan curve
  $\gamma:[0,1]\to\partial\Omega$; fix  four real numbers  $0<a<b<c<d<1$. Then
  there exist  a unique positive  real number $L$  and a unique  conformal map
  $\Phi:\Omega\to(0,L)×(0,1)$,  with natural  extension to  $\bar\Omega$, such
  that   $\Phi(\gamma(a))=i$,  $\Phi(\gamma(b))=0$,   $\Phi(\gamma(c))=L$  and
  $\Phi(\gamma(d))=L+i$.
  
  $L$  is called  \emph{extremal distance}  between $\partial_1=\gamma([a,b])$
  and  $\partial_2=\gamma([c,d])$  in  $\Omega$;  it  is  denoted  $\dd_\Omega
  (\partial_1, \partial_2)$.
\end{hogtheorem}

\begin{dem}
  For the proof of this result, and much more about conformal maps and related
  topics    (including   the   proofs    of   Propositions \ref{prp:extbound1}
  and \ref{prp:extbound2}), cf. \cite{ahlfors:confinv}.
\end{dem}

\emph{Examples:} The extremal distance between  both sides of length $a$ in an
$a×b$  rectangle  is  $b/a$. By  the  analyticity  of  the logarithm  in  $\CC
\setminus    (-\infty,0]$,    if    $\Omega=\{\rho   e^{i\theta}:    r<\rho<R,
0<\theta<\alpha\}$ with $0<r<R<\infty$  and $0<\alpha<2\pi$, then the extremal
distance  in $\Omega$  between both  circle arcs  is $\alpha^{-1}  \log(R/r)$. 
Finally, if  $L$ is  the extremal distance  in $\Omega$ between  two connected
parts  $\partial_1$ and  $\partial_2$ of  $\partial\Omega$, then  the extremal
distance between  the two components of  $\partial\Omega \setminus (\partial_1
\cup \partial_2)$ is $L^{-1}$.

\begin{prp}
  \label{prp:extbound1}
  Let  $\rho:\Omega\to[0,\infty)$   be  a  continuous   function,  and  denote
  $A_\rho(\Omega) = \int\!\!\!\int _\Omega  \rho^2$ and for any continuous arc
  $\gamma$ in  $\Omega$, $L_\rho(\gamma) = \int_\gamma \rho(z)  |\dd z|$ (this
  defines the  Riemannian metric associated  with $\rho$). Then we  have, thus
  giving  a justification to  the term  \emph{extremal length},  the following
  characterization of $\dd_\Omega$:
  $$\dd_\Omega(\partial_1,\partial_2) = \Sup_\rho \Inf_{\gamma:\partial_1
    \rightsquigarrow \partial_2} \frac {L_\rho(\gamma)^2} {A_\rho(\gamma)}$$  
  (where $\gamma:\partial_1  \rightsquigarrow \partial_2$ means  that $\gamma$
  is  a  continuous   path  in  $\Omega$  with  first   and  second  endpoints
  respectively in $\partial_1$ and $\partial_2$).
\end{prp}

In many cases, it is sufficient to apply this with a finite family of $\rho$'s
to  obtain  a  fairly good  lower  bound  for  $\dd_\Omega$ ---  usually  even
$\rho=1$, \ie\  taking the Euclidean  metric, is sufficient.  Another estimate
for $\dd_\Omega$ is the following:

\begin{prp}
  \label{prp:extbound2}
  Let  $L$ be  a  positive real  number  and $f_1$,  $f_2:[0,L]\to\RR$ be  two
  continuous  functions such  that for  all $t$  in $[0,L]$  we  have $f_1(t)<
  f_2(t)$.  Introduce  $\Omega=\{x+iy:   0<x<L,  f_1(x)<y<f_2(x)\}$,  and  let
  $\partial_1$  and   $\partial_2$  stand  for  the   vertical  components  of
  $\partial\Omega$. Then:
  $$\dd_\Omega(\partial_1,\partial_2) \geqslant \int_0^L \frac {\dd t}
  {f_2(t)-f_1(t)}.$$
  Moreover, if $f_1$ has a continuous derivative and $f_2=f_1+a$, then
  $$\dd_\Omega(\partial_1,\partial_2) \leqslant \frac L a
  \left[1+\|f_1'\|_\infty^2\right]\!.$$
\end{prp}

\subsection{Some topological tools}

In this section, all sets considered will be assumed non-empty.

\begin{df}
  If $A$  is a subset of  the set $\CC$ of  complex numbers (or  of any Banach
  space), note
  $$V_r(A)=\{x\in\CC: \dd(x,A) < r\} = A + \mathcal B(0,r);$$  
  if  $A$   and  $B$  are  two   bounded  subsets  of   $\CC$,  introduce  the
  \emph{Hausdorff distance} between $A$ and $B$ as
  $$\dd_H(A,B) = \Inf\{r: A\subset V_r(B), B\subset V_r(A)\}.$$
  
  It is easy  to see that $\dd_H$ is nonnegative  and satisfies the triangular
  inequality  (namely $\dd_H(A,B)\leqslant\dd_H(A,C)+\dd_H(C,B)$ for  any $A$,
  $B$, $C$);  moreover $\dd_H(A,B)=0$ if and  only if $\bar  A=\bar B$. Hence,
  $\dd_H$ defines  a metric topology on  the set of compact  subsets of $\CC$,
  known as the \emph{Hausdorff topology}.
\end{df}

We will need  the following standard property about  the Hausdorff topology on
the subsets of some fixed set, describing the compact case:

\begin{prp}
  \label{prp:compact}
  Let $K$ be a compact subset of  $\CC$. Then the set $\mathcal P_c(K)$ of all
  (non-empty) closed subsets of $K$, equipped with the topology induced by the
  Hausdorff distance, is compact.
\end{prp}

\emph{Remark:} It is  still true (and the proof is  basically the same, except
in obtaining the fact that $A$  is non-empty and closed) that for any complete
space $E$ the  set $\mathcal P_c(E)$ is complete. Moreover,  if $E$ is locally
compact, so is $\mathcal P_c(E)$.  However, it is generally not bounded, hence
not compact.

\bibliographystyle{siam}
\bibliography{Biblio}

\begin{thebibliography}{10}

\bibitem{ahlfors:confinv}
{\sc L.~V. Ahlfors}, {\em Conformal Invariants}, Topics in Geometric Function
  Theory, McGraw-Hill, New York, 1973.

\bibitem{bass:cutting}
{\sc R.~F. Bass and K.~Burdzy}, {\em Cutting {Brownian} Paths}, vol.~657 of
  Memoirs of the American Mathematical Society, AMS, 1999.

\bibitem{cranston:rw}
{\sc M.~C. Cranston and T.~S. Mountford}, {\em An extension of a result of
  {Burdzy} and {Lawler}}, Probab. Theory Relat. Fields, 89 (1991),
  pp.~487--502.

\bibitem{evans:cone}
{\sc S.~N. Evans}, {\em On the {Hausdorff} dimension of {Brownian} cone
  points}, Math. Proc. Camb. Phil. Soc., 98 (1985), pp.~343--353.

\bibitem{kahane:random}
{\sc J.-P. Kahane}, {\em Some Random Series of Functions}, vol.~5 of Cambridge
  Studies in Advanced Mathematics, Cambridge University Press, 2~ed., 1993.

\bibitem{kaufman:metrique}
{\sc R.~Kaufman}, {\em Une propriété métrique du mouvement brownien}, C.R.
  Acad. Sci., A 268 (1969), pp.~727--728.

\bibitem{lawler:frontier}
{\sc G.~F. Lawler}, {\em The dimension of the frontier of planar {Brownian}
  motion}, Elect. Comm. in Probab., 1 (1996), pp.~29--47.

\bibitem{lawler:hausdorff}
\leavevmode\vrule height 2pt depth -1.6pt width 23pt, {\em {Hausdorff}
  dimension of cut points for {Brownian} motion}, Electronic Journal of
  Probability, 1 (1996), pp.~1--20.

\bibitem{lawler:fractal}
\leavevmode\vrule height 2pt depth -1.6pt width 23pt, {\em Geometric and
  fractal properties of {Brownian} motion and random walk paths in two and
  three dimensions}, in Workshop in Random Walks, Budapest, 1998.

\bibitem{lawler:concavity}
\leavevmode\vrule height 2pt depth -1.6pt width 23pt, {\em Strict concavity of
  the intersection exponent for {Brownian} motion in two and three dimensions},
  Mathematical Physics Electronic Journal, 4 (1998).

\bibitem{lawler:walk}
{\sc G.~F. Lawler and E.~E. Puckette}, {\em The intersection exponent for
  simple random walks}, Comb. Prob. Comp., 9 (2000), pp.~441--464.

\bibitem{werner:analyticity}
{\sc G.~F. Lawler, O.~Schramm, and W.~Werner}, {\em Analyticity of intersection
  exponents for planar {Brownian} motion}.
\newblock to appear, 2000.

\bibitem{werner:value}
\leavevmode\vrule height 2pt depth -1.6pt width 23pt, {\em Values of {Brownian}
  intersection exponents {I}: Half-plane exponents}, Acta Mathematica,  (2000).
\newblock To appear.

\bibitem{werner:value2}
\leavevmode\vrule height 2pt depth -1.6pt width 23pt, {\em Values of {Brownian}
  intersection exponents {II}: Plane exponents}, Acta Mathematica,  (2000).
\newblock To appear.

\bibitem{werner:value3}
\leavevmode\vrule height 2pt depth -1.6pt width 23pt, {\em Values of {Brownian}
  intersection exponents {III}: Two-sided exponents}, Annales Inst. Henri
  Poincar\'e,  (2000).
\newblock to appear.

\bibitem{werner:4/3}
\leavevmode\vrule height 2pt depth -1.6pt width 23pt, {\em The dimension of the
  {Brownian} frontier is $4/3$.}, Math. Res. Lett., 8 (2001), pp.~13--24.

\bibitem{lawler:intersection}
{\sc G.~F. Lawler and W.~Werner}, {\em Intersection exponents for planar
  {Brownian} motion}, Ann. Probab.,  (1999), pp.~1601--1642.

\bibitem{legall:sflour}
{\sc J.-F. Le~Gall}, {\em Some properties of planar {Brownian} motion}, in
  École d'été de Probabilités de Saint-Flour XX--1990, vol.~1527 of Lecture
  Notes in Mathematics, Springer, 1992.

\bibitem{schramm:UST}
{\sc O.~Schramm}, {\em Scaling limits of loop-erased random walks and uniform
  spanning trees}, Israel Journal of Mathematics, 118 (2000), pp.~221--288.

\bibitem{smirnov:perco}
{\sc S.~Smirnov}, {\em Critical percolation in the plane}, Preprint,  (2001).

\bibitem{werner:disconnection}
{\sc W.~Werner}, {\em On {Brownian} disconnection exponents}, Bernoulli, 1
  (1995), pp.~371--380.

\bibitem{werner:bounds}
\leavevmode\vrule height 2pt depth -1.6pt width 23pt, {\em Bounds for
  disconnection exponents}, Elect. Comm. in Probab., 1 (1996), pp.~19--28.

\bibitem{werner:barcelone}
\leavevmode\vrule height 2pt depth -1.6pt width 23pt, {\em Critical exponents,
  conformal invariance and planar {Brownian} motion}, in Proceedings of the 3rd
  {European} {Mathematical} {Congress}, Birkhäuser, 2000.
\newblock To appear.

\end{thebibliography}

\end{document}